\documentclass[12pt]{amsart}
\usepackage[colorlinks=true,citecolor=black,linkcolor=black,urlcolor=blue]{hyperref}
\usepackage{amsmath}
\usepackage[english,  activeacute]{babel}
\usepackage[utf8]{inputenc}
\usepackage{amssymb}
\usepackage{amsthm}
\usepackage{graphics,graphicx}
\usepackage{enumerate}
\usepackage{ytableau}
\usepackage{array}
\usepackage{bm}
\usepackage{cite}
\usepackage{a4wide}
\usepackage{color, url}
\usepackage{float}
\usepackage{comment}
\usepackage{xcolor}
\usepackage{color, url}
\usepackage{float}
\usepackage{accents}
\usepackage{pgf, tikz}
\usepackage{multirow}
\usepackage{tikz}

\usepackage{booktabs}
\usepackage{tabularx}
\usepackage{geometry}
\theoremstyle{plain}
\newtheorem{theorem}{Theorem}[section]

\newtheorem{lemma}[theorem]{Lemma}
\newtheorem{remark}[theorem]{Remark}

\theoremstyle{definition}

\newcommand{\Q}{{\mathbb Q}}

\newcommand{\lth}[1]{\texttt{length}#1}
\newcommand{\last}[1]{\texttt{last}#1}
\newcommand{\bck}[1]{\texttt{black}#1}
\newcommand{\ver}[1]{\texttt{verblack}#1}
\newcommand{\white}[1]{\texttt{verwhite}#1}

\usepackage[usestackEOL]{stackengine}[2013-10-15]
\def\x{\hspace{3ex}}    
\def\y{\hspace{2.45ex}}  
\stackMath

\newcommand{\jluc}[1]{\mbox{}{\sf\color{blue}[jluc: #1]}\marginpar{\color{blue}\Large$*$}} 
\title{Black Cell Capacity in Catalan polyominoes}

\author[J.-L. Baril]{Jean-Luc Baril}
\address{LIB, Universit\'e Bourgogne Europe,
  B.P. 47 870, 21078 Dijon Cedex France}
\email{barjl@u-bourgogne.fr}

\author[S. Fried]{Sela Fried}
\address{ Department of Computer Science, Israel Academic College, 52275 Ramat Gan, Israel}
\email{friedsela@gmail.com}

\author[N. Hassler]{Nathana\"el Hassler}
\address{LIB, Universit\'e Bourgogne Europe,
  B.P. 47 870, 21078 Dijon Cedex France}
\email{nathanael.hassler@ens-rennes.fr}

\author[J.L. Ram\'irez]{Jos\'e L. Ram\'irez}
\address{Departamento de Matem\'aticas,  Universidad Nacional de Colombia,  Bogot\'a, Colombia}
\email{jlramirezr@unal.edu.co}

\date{\today}

\begin{document}

\newcommand{\nadji}[1]{\mbox{}{\sf\color{green}[Ram\'{\i}rez: #1]}\marginpar{\color{green}\Large$*$}} 

\tikzstyle{carre}=[draw, minimum size=1cm]
\newcommand{\case}[2]{
    \pgfmathparse{mod(#1+#2,2)==0 ? "black" : "white"}
    \edef\couleur{\pgfmathresult}
    \filldraw[fill=\couleur] (#1,#2) rectangle ++(1,1);
}
\tikzstyle{carre}=[draw, minimum size=1cm]
\newcommand{\casee}[2]{
    \pgfmathparse{mod(#1+#2,2)==0 ? "white" : "black"}
    \edef\couleur{\pgfmathresult}
    \filldraw[fill=\couleur] (#1,#2) rectangle ++(1,1);
}

\begin{abstract} A Catalan word is a sequence $w_1w_2\cdots w_n$ of nonnegative integers such that $w_1=0$ and $w_{i}\leq w_{i-1}+1$ for $2\leq i\leq n$. Given a Catalan word, we construct a column-convex polyomino (or \emph{bargraph}) by placing, at position $i$, a column of height $w_i + 1$, with all columns aligned along their bottom edges. On these Catalan polyominoes we define the black cell capacity by coloring the cells in a chessboard pattern and we count the number of black cells in the polyomino.  We study the distribution of the black cell capacity over Catalan polyominoes and derive generating functions that encode this statistic.
\end{abstract}

\maketitle

\section{Introduction}

A \emph{Catalan word} of length $n \geq 1$ is a sequence $w = w_1 w_2 \cdots w_n$ of nonnegative integers with $w_1 = 0$ and $w_i \leq w_{i-1} + 1$ for $i = 2,\ldots,n$. For $n=0$, the unique Catalan word is the empty word $\epsilon$.  Let $\mathcal{C}_n$ denote the set of all Catalan words of length $n$. It is well known that $|\mathcal{C}_n|=C_n$, where $C_n$ is the $n$th \emph{Catalan number},
\[
C_n = \frac{1}{n+1}\binom{2n}{n},
\]
see, for instance,  \cite[Exercise~80]{stan}. For enumerations of Catalan words with respect to various statistics, see   \cite{Baril3,BGR,Baril2,relation,CallManRam,AlejaRam}.

Catalan words are closely related to certain lattice paths in the first quadrant. More precisely, a \emph{Dyck path} of semilength $n$ is a lattice path in  the first quarter plane that starts at $(0,0)$, ends  at $(2n,0)$, and uses only up-steps $U=(1,1)$ and down-steps $D=(1,-1)$. Such a path can be encoded by a word in $\mathcal{C}_n$ by recording, from left to right, the $y$-coordinates of the initial vertices of its up-steps.  Background on lattice paths and Dyck paths can be found in \cite{Banfla,Deu}.

Catalan words also admit a natural geometric representation. Given a word $w = w_1 \cdots w_n$, one obtains a column-convex polyomino (or \emph{bargraph}) by placing, at position $i$, a column of height $w_i+1$, with all columns aligned along their bottom edges.  The resulting object is called a \emph{Catalan polyomino}. Let $\mathbf{C}_n$ be the set of these polyominoes with $n$ columns, and let $\mathbf{C}=\cup_{n\geq 1}\mathbf{C}_n$.  We refer to \cite{Book1} for a historical review of polyominoes, and to \cite{BLE3,ManSha2,ManRam} for definitions and enumerative methods related to polyominoes. See also \cite{BlKn} for work on Catalan polyominoes.

For a Catalan polyomino $P\in\mathbf{C}$, we denote by $\lth(P)$ the number of columns of $P$, which is called the {\it length} of $P$. The number of cells in its last column will be denoted $\last(P)$. Following  \cite{Chen,fried}, we color the cells of $P$ in a chessboard pattern, with the southwestern cell colored black, and we study the distribution of the {\it black cell capacity} $\bck(P)$, that is, the number of black cells contained in $P$.  

\begin{figure}[ht!]
\centering
\begin{tikzpicture}[scale=0.5]-
\filldraw[fill=black] (0,0) rectangle ++(1,1);
\foreach \x in {1,2,3,4,5,6,7,8,9,10,11,12} {
    \case{\x}{0}
}
\foreach \y in {1} {
    \case{2}{\y}
}
\foreach \y in {1,2} {
    \case{3}{\y}
}
\foreach \y in {1} {
    \case{5}{\y}
}
\foreach \y in {1,2} {
    \case{6}{\y}
}
\foreach \y in {1,2,3} {
    \case{7}{\y}
}
\foreach \y in {1} {
    \case{8}{\y}
}
\foreach \y in {1} {
    \case{10}{\y}
}
\foreach \y in {1} {
    \case{11}{\y}
}
\end{tikzpicture}\qquad\quad
\begin{tikzpicture}[scale=0.5]-
\filldraw[fill=black] (0,0) rectangle ++(1,1);
\foreach \x in {1,2,3,4,5,6,7,8,9,10,11,12} {
    \case{\x}{0}
}
\foreach \y in {1} {
    \casee{2}{\y}
}
\casee{3}{1}
\case{3}{2}
\foreach \y in {1} {
    \casee{5}{\y}
}
    \casee{6}{1}
    \case{6}{2}
     \casee{7}{1}
    \case{7}{2}
    \casee{7}{3}

 \casee{8}{1}
  
\foreach \y in {1} {
    \casee{10}{\y}
}
\foreach \y in {1} {
    \casee{11}{\y}
}
\end{tikzpicture}
\caption{The polyomino $P$ associated with the Catalan Word $w=0012012310110$. We have $\lth(P)=13$, $\last(P)=1=w_{13}+1$. The black cell capacity $\bck(P)$  equals to 13,  the vertical black cell capacity $\ver(P)$ equals to 12, and the vertical white cell capacity $\white(P)$ equals to 13. } \label{polyo}
\end{figure}
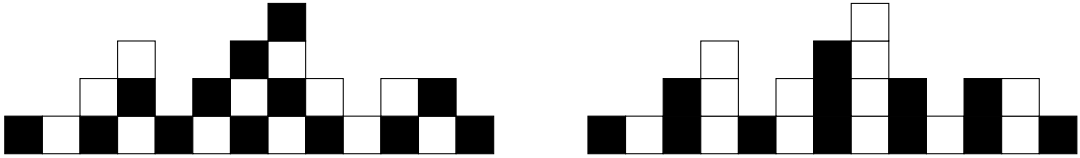

We also define the {\it vertical black cell capacity} of a polyomino $P$, denoted $\ver(P)$, as the  total number of cells in columns of odd index (with the first column having index 1). Similarly, the {\it vertical white cell capacity}, denoted $\white(P)$, is the total number of cells in columns of even index. Figure~\ref{polyo} shows the Catalan polyomino of length $13$ associated with the word \texttt{0012012310110}, where $\bck(P)=13$, $\ver(P)=12$, and $\white(P)=13$.

\noindent{\bf Black cell capacity versus vertical black/white cell capacity.}   Let $P$ be a Catalan polyomino of length $n\geq 1$, colored in a chessboard pattern. Then $P$ decomposes into $n$ \emph{northeast diagonals} of cells, each consisting entirely of cells of the same color, such that the $i$th diagonal and the $i$th column of $P$ share the same bottom cell. Let $D_i$ be the $i$th diagonal for $1\leq i\leq n$, and let $a_i$ be the number of cells in $D_i$. From $P$ we construct a polyomino $Q$ as follows: the columns of $Q$ are aligned along a common baseline, and the $i$th column of $Q$ has exactly $a_{n-i+1}$ cells. This construction defines a bijection from $\mathbf{C}_n$ to itself; denote it by $f$. By construction, $\bck(P)=\ver(f(P))$ if $n$ is odd, and 
$\bck(P)=\white(f(P))$ if $n$ is even. See Figures~\ref{bijex1} and \ref{bijex2} for illustrations of $f$ in the cases $n$ odd and $n$ even.

\begin{figure}[ht!]
\centering
\begin{tikzpicture}[scale=0.5]-
\filldraw[fill=black] (0,0) rectangle ++(1,1);
\foreach \x in {1,2,3,4,5,6} {
    \case{\x}{0}
}
\foreach \y in {1} {
    \case{1}{\y}
}
\foreach \y in {1,2} {
    \case{2}{\y}
}
\foreach \y in {1} {
    \case{3}{\y}
}
\foreach \y in {1} {
    \case{5}{\y}
}
\end{tikzpicture}\quad $\longrightarrow$ \quad
\begin{tikzpicture}[scale=0.5]-
\filldraw[fill=black] (0,0) rectangle ++(1,1);
\foreach \x in {1,2,3,4,5,6} {
    \case{\x}{0}
}
\casee{2}{1}
\casee{4}{1}
\casee{5}{1}
\casee{6}{1}\case{6}{2}
\end{tikzpicture}
\caption{ The image by $f$ of the Catalan polyomino $P=1232121$ is $f(P)=1121223$. The length of $P$ is odd and thus, we have $\bck(P)=\ver(f(P))=8$.} \label{bijex1}
\end{figure}
\begin{figure}[H]
\centering
\begin{tikzpicture}[scale=0.5]-
\filldraw[fill=black] (0,0) rectangle ++(1,1);
\foreach \x in {1,2,3,4,5,6,7} {
    \case{\x}{0}
}
\foreach \y in {1} {
    \case{1}{\y}
}
\foreach \y in {1,2} {
    \case{2}{\y}
}
\foreach \y in {1} {
    \case{3}{\y}
}
\foreach \y in {1} {
    \case{5}{\y}
\case{7}{1}
}
\end{tikzpicture}\quad $\longrightarrow$ \quad
\begin{tikzpicture}[scale=0.5]-
\filldraw[fill=black] (0,0) rectangle ++(1,1);
\foreach \x in {1,2,3,4,5,6,7} {
    \case{\x}{0}
}
\casee{1}{1}
\casee{3}{1}
\casee{5}{1}
\casee{6}{1}
\casee{7}{1}\case{7}{2}
\end{tikzpicture}
\caption{ The image by $f$ of the Catalan polyomino $P=12321212$ is $f(P)=12121223$. The length of $P$ is  even and thus, we have  $\bck(P)=\white(f(P))=9$.} \label{bijex2}
\end{figure}

Let us formalize the bijection stated above for later use.

\begin{theorem} There is a bijection $f$ on Catalan words  of odd length (resp. even length) that transports the black cell capacity into the vertical black cell capacity (resp. vertical white capacity).
\label{th11}
\end{theorem}

\noindent{\bf Motivation and Outline of the paper.} 
In \cite{Chen}, the authors prove that the total number of black cells under all Dyck paths of semilength $n$ is $4^{n-1}-{2n \choose n-1}$. However, they do not provide any results on the distribution of the black cell capacity.

In this paper, we study the distribution of black cell statistics on Catalan polyominoes. In Section 2, we derive a system of four functional equations for the generating functions $F_{ab}(x,u,q)$, $a,b\in\{0,1\}$. These series count Catalan polyominoes of length congruent to $a$ modulo 2  (marked by the variable $x$), ending with a column containing  $b$ cells modulo 2 (marked by the variable $u$), and classified by the black cell capacity (marked by the variable $q$). We express the system as a matrix equation, which yields a formal solution that can be used, together with computer algebra, to compute the initial terms of the corresponding series expansions.

In Section 3, we take the continued fraction approach. The solution can be expressed as a $2 \times 2$ matrix continued fraction, whose $n$th convergent provides the initial terms of the series expansion.

In Section 4, we  investigate the bistatistic 
$(r,s)$, where $r$ counts the cells in odd-indexed columns and  $s$ counts the cells in even-indexed columns. This approach allows us to derive a closed form for the generating function of Catalan polyominoes with respect to their length, the bistatistic 
$(r,s)$, and the size of the last column. From this, we deduce  a closed form for the generating function encoding the distribution of the black cell capacity.

In Section 5, we formulate a functional equation (see Theorem~\ref{thm51}) for the distribution of the black cell capacity and we solve it by giving a closed form of the generating function.

\section{Length, black cell capacity, and last value} 

In this section we count Catalan polyominoes by  length, black cell capacity, and the size of the last column.

Let $n,i, k\in\mathbb{N}$ with $n,i,k\geq 1$. We denote by $\mathbf{C}_{n,i,k}$ the set of Catalan polyominoes of length $n$ whose  $n$th (last) column has $i$ cells, and whose black cell capacity equals $k$. Set
\[
F_{n,i,k}:=|\mathbf{C}_{n,i,k}|,
\]
and define the trivariate ordinary generating function
$$F(x,u,q)\;=\;\sum_{n,i,k\geq 1} F_{n,i,k}x^{n}u^{i}q^{k}
=\;\sum_{P\in\mathbf{C}}x^{\lth(P)}u^{\last(P)}q^{\bck(P)}.$$
Here $x$ marks the length (number of columns), $u$ marks the number of cells in the last column, and $q$ marks the black cell capacity. For $a,b\in\{0,1\}$, we denote by $F_{ab}(x,u,q)$  the generating function of the subclass of Catalan polyominoes $P$ such that $\lth(P)=a \mod{2}$ and $\last(P)=b \mod 2$.
In particular,
$$F_{00}(x,u,q)=\sum_{n,i,k\geq 1} F_{2n,2i,k}\,x^{2n}u^{2i}q^{k}, \quad F_{10}(x,u,q)=\sum_{n,i,k\geq 1} F_{2n-1,2i,k}\,x^{2n-1}u^{2i}q^{k},$$
$$F_{01}(x,u,q)=\sum_{n,i,k\geq 1} F_{2n,2i-1,k}\,x^{2n}u^{2i-1}q^{k}, \quad F_{11}(x,u,q)=\sum_{n,i,k\geq 1} F_{2n-1,2i-1,k}\,x^{2n-1}u^{2i-1}q^{k}.$$
Clearly, 
$$F(x,u,q)=F_{00}(x,u,q)+F_{01}(x,u,q)+F_{10}(x,u,q)+F_{11}(x,u,q).$$

The next theorem provides a system of functional equations for the generating functions $F_{00}(x,u,q)$, $F_{10}(x,u,q)$, $F_{01}(x,u,q)$,  and $F_{11}(x,u,q)$.

\begin{theorem}\label{pol:3}
The generating functions $F_{ab}(x,u,q)$ for $a,b\in\{0,1\}$ satisfy the following system of functional equations: 
\[
\begin{aligned}
F_{00}(x,u,q)&=\frac{xu^2q}{qu^2-1}\left(F_{10}(x,\sqrt{q}u,q)-F_{10}(x,1,q) +\sqrt{q}uF_{11}(x,\sqrt{q}u,q)-F_{11}(x,1,q)   \right),\\
F_{10}(x,u,q)&=\frac{xu^2q}{qu^2-1}\left(F_{00}(x,\sqrt{q}u,q)-F_{00}(x,1,q)+\sqrt{q}uF_{01}(x,\sqrt{q}u,q)-F_{01}(x,1,q)\right),\\
F_{01}(x,u,q)&=\frac{xu}{qu^2-1}\left(u^2qF_{10}(x,\sqrt{q}u,q)-F_{10}(x,1,q)+\sqrt{q}uF_{11}(x,\sqrt{q}u,q)-F_{11}(x,1,q)\right),\\
F_{11}(x,u,q)&=xuq+\frac{xuq}{qu^2-1}\left(u^2qF_{00}(x,\sqrt{q}u,q)-F_{00}(x,1,q)\right)+\\
&\hskip6cm\frac{xuq}{qu^2-1}\left(\sqrt{q}uF_{01}(x,\sqrt{q}u,q)-F_{01}(x,1,q)\right).
\end{aligned}
\]
\end{theorem}

\begin{proof}
We distinguish four types of Catalan polyominoes $P$ according to the parities of the number \lth(P) of columns  and the number \last(P) of cells in the last column. Since the four cases are similar, we only give the proof for the case in which $\lth(P)$ and $\last(P)$ are odd.

\smallskip
Let us assume that $\lth(P)$ and $\last(P)$ are odd. If $P$ has a single column, then it consists of one cell, and its contribution to the generating function is $xuq$. Otherwise, assume that $P$ has at least two columns. Let $Q$ be the polyomino obtained from $P$ by deleting its last column. Then 
$\lth(Q)=\lth(P)-1=0\mod{2}$ and the Catalan structure of $P$ implies the condition $\last(Q)+1\geq \last(P)$.

\begin{figure}[ht!]
    \centering
\begin{tikzpicture}[scale=0.5]-
\filldraw[fill=black] (0,0) rectangle ++(1,1);
\foreach \x in {1,2,3} {
    \case{\x}{0}
}
\draw[dashed] (4.2,0.5) -- (5.8,0.5);
\casee{6}{0}
\casee{7}{0}
\foreach \y in {1,2,3} {
    \casee{6}{\y}
}
\foreach \y in {1,2,3,4} {
    \casee{7}{\y}
}
\draw[thick]
  plot[domain=1:4, smooth]
  (\x, {1 + (\x-1) + sin((\x-1)*200)});
\draw[dotted] (4,0)--(4,3.3);
\end{tikzpicture}\qquad \qquad\qquad \qquad
\begin{tikzpicture}[scale=0.5]-
\filldraw[fill=black] (0,0) rectangle ++(1,1);
\foreach \x in {1,2,3} {
    \case{\x}{0}
}
\draw[dashed] (4.2,0.5) -- (5.8,0.5);
\casee{7}{0}
\casee{6}{0}
\foreach \y in {1,2} {
    \casee{6}{\y}
}
\foreach \y in {1,2} {
    \casee{7}{\y}
}
\draw[thick]
  plot[domain=1:4, smooth]
  (\x, {1 + (\x-1) + sin((\x-1)*200)});
\draw[dotted] (4,0)--(4,3.3);
\end{tikzpicture}
\caption{Illustration of Case 1: $\lth(P)$ and $\last(P)$ are odd. The left part shows a polyomino $P$ where $Q$ satisfies $\last(Q)=0\mod 2$, while the right part is for $\last(Q)=1\mod 2$. }
\label{case1}
\end{figure}
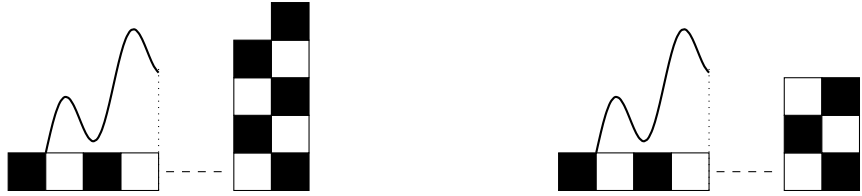

If $\last(Q)$ is even (see Figure~\ref{case1}), then the contribution for these polyominoes is 
\begin{align*}
    A_1:=& x\sum_{n,i,k\geq 1} F_{2n,2i,k}\,x^{2n}q^k\left(uq+ u^3q^2+u^5q^3+\cdots + u^{2i+1}q^{i+1}\right)\\
    =& xuq\sum_{n,i,k\geq 1} F_{2n,2i,k}\,x^{2n}q^k\left(\frac{u^{2i+2}q^{i+1}-1}{qu^2-1}\right)\\
    =& \frac{xuq}{qu^2-1}\left(u^2qF_{00}(x,\sqrt{q}u,q)-F_{00}(x,1,q)\right).
\end{align*}
If $\last(Q)$ is odd (see Figure~\ref{case1}), then the contribution for these polyominoes is 
\begin{align*}
B_1:=& x\sum_{n,i,k\geq 1} F_{2n,2i-1,k}\,x^{2n}q^k \left(uq+ u^3q^2+u^5q^3+\cdots + u^{2i-1}q^i\right)\\
=&xuq \sum_{n,i,k\geq 1} F_{2n,2i-1,k}\,x^{2n}q^k \left(\frac{u^{2i}q^i-1}{qu^2-1}\right)\\
    =&\frac{xuq}{qu^2-1}\left(\sqrt{q}uF_{01}(x,\sqrt{q}u,q)-F_{01}(x,1,q)\right).
\end{align*}
For this case, the contribution is therefore
\begin{align*} F_{11}(x,u,q)=&xuq+A_1+B_1\\
&xuq+\frac{xuq}{qu^2-1}\left(u^2qF_{00}(x,\sqrt{q}u,q)-F_{00}(x,1,q)\right.+\\
&\hskip4cm\left.\sqrt{q}uF_{01}(x,\sqrt{q}u,q)-F_{01}(x,1,q)\right).
\end{align*}

\end{proof}

Now, let us define the generating functions $G_{ab}(x,u,q)$ for $a,b\in\{0,1\}$ as follows:
\begin{align*}
  G_{00}(x,u,q)=&F_{00}(x,\sqrt{u},q), \qquad \qquad
  G_{01}(x,u,q)=\sqrt{u}F_{01}(x,\sqrt{u},q)\\
  G_{10}(x,u,q)=&F_{10}(x,\sqrt{u},q), \qquad \qquad
  G_{11}(x,u,q)=\sqrt{u}F_{11}(x,\sqrt{u},q).
\end{align*}
We also set \[\mathbf{G}(x,u,q)=
\begin{pmatrix} G_{00}(x,u,q)\\
                G_{01}(x,u,q)\\
                G_{10}(x,u,q)\\
                G_{11}(x,u,q)
\end{pmatrix}.\] 

Considering these new generating functions, it is straightforward to see that Theorem~2.1 can also be written as follows:

\begin{theorem}\label{pol:3}
The generating functions $G_{ab}(x,u,q)$ for $a,b\in\{0,1\}$ satisfy the following system of functional equations: 
\[
\begin{aligned}
G_{00}(x,u,q)&=\frac{xuq}{qu-1}\left(G_{10}(x,qu,q)-G_{10}(x,1,q)+G_{11}(x,qu,q)-G_{11}(x,1,q)\right),\\
G_{01}(x,u,q)&=\frac{xu}{qu-1}\left(uqG_{10}(x,qu,q)-G_{10}(x,1,q)+G_{11}(x,qu,q)-G_{11}(x,1,q)\right),\\
G_{10}(x,u,q)&=\frac{xuq}{qu-1}\left(G_{00}(x,qu,q)-G_{00}(x,1,q)+G_{01}(x,qu,q)-G_{01}(x,1,q)\right),\\
G_{11}(x,u,q)&=xuq+\frac{xuq}{qu-1}\left(uqG_{00}(x,qu,q)-G_{00}(x,1,q)+G_{01}(x,qu,q)-G_{01}(x,1,q)\right),
\end{aligned}
\]
 which is equivalent to the following matrix equation:
 \begin{align}
 \mathbf{G}(x,u,q)=\mathbf{M}(x,u,q)\cdot \mathbf{G}(x,qu,q)-\mathbf{N}(x,u,q)\cdot  \mathbf{G}(x,1,q)+\mathbf{B}(x,u,q).
 \end{align} 
 where

\[\mathbf{M}(x,u,q)=
\begin{pmatrix}
0& 0 & \frac{xuq}{qu-1} & \frac{xuq}{qu-1} \\ 
 0& 0 & \frac{xu^2q}{qu-1} & \frac{xu}{qu-1}\\  
 \frac{xuq}{qu-1} & \frac{xuq}{qu-1} & 0 & 0 \\ 
\frac{xu^2q^2}{qu-1} & \frac{xuq}{qu-1} & 0 & 0  
\end{pmatrix},
 \quad \mathbf{N}(x,u,q)=
\begin{pmatrix}
0& 0 & \frac{xuq}{qu-1} & \frac{xuq}{qu-1} \\ 
 0& 0 & \frac{xu}{qu-1} & \frac{xu}{qu-1}\\  
 \frac{xuq}{qu-1} & \frac{xuq}{qu-1} & 0 & 0 \\ 
\frac{xuq}{qu-1} & \frac{xuq}{qu-1} & 0 & 0  
\end{pmatrix},
\]

\[
\mathbf{B}(x,u,q)=
\begin{pmatrix} 0\\
                0\\
                0\\
                xuq\\
\end{pmatrix}.
\]
\end{theorem}

By iterating Equation (1), we obtain a formal expression for the vector $\mathbf{G}(x,u,q)$.

\begin{theorem} We have
 \[
\mathbf{G}(x,u,q) = \sum_{k=0}^{\infty} \mathbf{P}_k(x,u,q) \cdot \Big( \mathbf{B}(x,q^k u,q) - \mathbf{N}(x,q^k u,q) \, \mathbf{G}(x,1,q) \Big),
\]
with 
\[
\mathbf{P}_k(x,u,q) := \prod_{j=0}^{k-1} \mathbf{M}(x,q^j u,q), \quad \mathbf{P}_{-1}(x,u,q)=\mathbf{P}_0(x,u,q)= \mathbf{I},
\]
where $\mathbf{I}$ is the identity matrix of size 4.
\end{theorem}

Substituting $u=1$ in the expression of $\mathbf{G}(x,u,q)$, we obtain:
\[
\mathbf{G}(x,1,q) = \sum_{k=0}^{\infty} \mathbf{P}_k(x,1,q) \cdot\Big( \mathbf{B}(x,q^k,q) - \mathbf{N}(x,q^k,q) \, \mathbf{G}(x,1,q) \Big).
\]
Isolating $\mathbf{G}(x,1,q)$ we obtain
\[
\left( \mathbf{I} + \sum_{k=0}^{\infty} \mathbf{P}_k(x,1,q) \mathbf{N}(x,q^k,q) \right) \cdot \mathbf{G}(x,1,q) = \sum_{k=0}^{\infty} \mathbf{P}_k(x,1,q) \mathbf{B}(x,q^k,q),
\]
and thus
\[
 \mathbf{G}(x,1,q) = \left( \mathbf{I} + \sum_{k=0}^{\infty} \mathbf{P}_k(x,1,q) \mathbf{N}(x,q^k,q) \right)^{-1} \cdot\sum_{k=0}^{\infty} \mathbf{P}_k(x,1,q) \mathbf{B}(x,q^k,q).
\]
Then, we deduce the following theorem.

\begin{theorem}
\[
 \mathbf{G}(x,1,q)=\left( \sum_{n=0}^\infty (-1)^n\mathbf{S}(x,q)^n\right)\cdot\sum_{k=0}^{\infty} \mathbf{P}_k(x,1,q) \mathbf{B}(x,q^k,q),
 \]
where
\[\mathbf{S}(x,q)=\sum_{k=0}^{\infty} \mathbf{P}_k(x,1,q) \mathbf{N}(x,q^k,q).\]
\end{theorem}
Notice that if we denote by $\mathbf{F}(x,u,q)$ the vector \[
\mathbf{F}(x,u,q)=
\begin{pmatrix} F_{00}(x,u,q)\\
                F_{01}(x,u,q)\\
                F_{10}(x,u,q)\\
                F_{11}(x,u,q)
\end{pmatrix},
\]
then we have $\mathbf{F}(x,1,q)=\mathbf{G}(x,1,q)$.

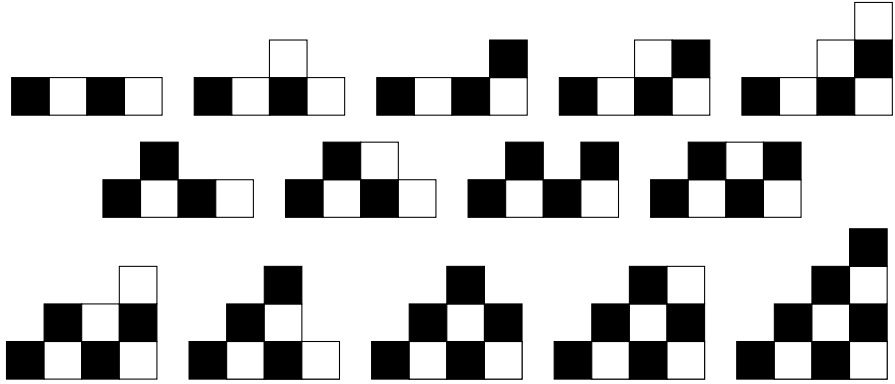
\begin{figure}[ht!]
\centering
\begin{tikzpicture}[scale=0.5]-
\filldraw[fill=black] (0,0) rectangle ++(1,1);
\foreach \x in {1,2,3} {
    \case{\x}{0}
}
\end{tikzpicture}\quad 
\begin{tikzpicture}[scale=0.5]-
\filldraw[fill=black] (0,0) rectangle ++(1,1);
\foreach \x in {1,2,3} {
    \case{\x}{0}
}    \case{2}{1}
\end{tikzpicture}\quad 
\begin{tikzpicture}[scale=0.5]-
\filldraw[fill=black] (0,0) rectangle ++(1,1);
\foreach \x in {1,2,3} {
    \case{\x}{0}
}
    \case{3}{1}
\end{tikzpicture}\quad 
\begin{tikzpicture}[scale=0.5]-
\filldraw[fill=black] (0,0) rectangle ++(1,1);
\foreach \x in {1,2,3} {
    \case{\x}{0}
}
\case{2}{1}
    \case{3}{1}
\end{tikzpicture}\quad 
\begin{tikzpicture}[scale=0.5]-
\filldraw[fill=black] (0,0) rectangle ++(1,1);
\foreach \x in {1,2,3} {
    \case{\x}{0}
}
    \case{2}{1}
\foreach \y in {1,2} {
    \case{3}{\y}
}
\end{tikzpicture}\\[0.75em]
\begin{tikzpicture}[scale=0.5]-
\filldraw[fill=black] (0,0) rectangle ++(1,1);
\foreach \x in {1,2,3} {
    \case{\x}{0}
}
    \case{1}{1}
\end{tikzpicture}\quad 
\begin{tikzpicture}[scale=0.5]-
\filldraw[fill=black] (0,0) rectangle ++(1,1);
\foreach \x in {1,2,3} {
    \case{\x}{0}
}
    \case{1}{1}
    \case{2}{1}
\end{tikzpicture}\quad 
\begin{tikzpicture}[scale=0.5]-
\filldraw[fill=black] (0,0) rectangle ++(1,1);
\foreach \x in {1,2,3} {
    \case{\x}{0}
}
    \case{1}{1}
    \case{3}{1}
\end{tikzpicture}\quad 
\begin{tikzpicture}[scale=0.5]-
\filldraw[fill=black] (0,0) rectangle ++(1,1);
\foreach \x in {1,2,3} {
    \case{\x}{0}
}
    \case{1}{1}
    \case{2}{1}
    \case{3}{1}
\end{tikzpicture}\\[0.25em]
\begin{tikzpicture}[scale=0.5]-
\filldraw[fill=black] (0,0) rectangle ++(1,1);
\foreach \x in {1,2,3} {
    \case{\x}{0}
}
    \case{1}{1}
    \case{2}{1}
\foreach \y in {1,2} {
    \case{3}{\y}
}
\end{tikzpicture}\quad
\begin{tikzpicture}[scale=0.5]-
\filldraw[fill=black] (0,0) rectangle ++(1,1);
\foreach \x in {1,2,3} {
    \case{\x}{0}
}
\foreach \y in {0,1} {
    \case{1}{\y}
}
\foreach \y in {0,1,2} {
    \case{2}{\y}
}
\foreach \y in {0} {
    \case{3}{\y}
}
\end{tikzpicture}\quad 
\begin{tikzpicture}[scale=0.5]-
\filldraw[fill=black] (0,0) rectangle ++(1,1);
\foreach \x in {1,2,3} {
    \case{\x}{0}
}
\foreach \y in {0,1} {
    \case{1}{\y}
}
\foreach \y in {0,1,2} {
    \case{2}{\y}
}
\foreach \y in {0,1} {
    \case{3}{\y}
}
\end{tikzpicture}\quad 
\begin{tikzpicture}[scale=0.5]-
\filldraw[fill=black] (0,0) rectangle ++(1,1);
\foreach \x in {1,2,3} {
    \case{\x}{0}
}
\foreach \y in {0,1} {
    \case{1}{\y}
}
\foreach \y in {0,1,2} {
    \case{2}{\y}
}
\foreach \y in {0,1,2} {
    \case{3}{\y}
}
\end{tikzpicture}\quad 
\begin{tikzpicture}[scale=0.5]-
\filldraw[fill=black] (0,0) rectangle ++(1,1);
\foreach \x in {1,2,3} {
    \case{\x}{0}
}
\foreach \y in {0,1} {
    \case{1}{\y}
}
\foreach \y in {0,1,2} {
    \case{2}{\y}
}
\foreach \y in {0,1,2, 3} {
    \case{3}{\y}
}
\end{tikzpicture}
\caption{The 14 Catalan polyominoes of length 4. There are 2 (resp 5, 4,2,1) polyominoes with black cell capacity 2 (resp. 3,4,5,6).} 
\label{polyo14}
\end{figure}

The first terms of the series expansion of the generating function  \begin{align*}F(x,1,q)=&G(x,1,q)=\begin{pmatrix}
1 & 1 & 1 &  1
\end{pmatrix}\cdot \mathbf{G}(x,1,q)\\
=&F_{00}(x,1,q)+F_{01}(x,1,q)+F_{10}(x,1,q)+F_{11}(x,1,q)\end{align*} are  
\begin{multline*}
xq+(q^2+q)x^2+(q^4+2q^3+2q^2)x^3+\bm{(q^6+2q^5+4q^4+5q^3+2q^2)x^4}+\\
(q^9+2q^8+5q^7+8q^6+12q^5+10q^4+4q^3)x^5+\\
(q^{12}+2q^{11}+5q^{10}+9q^9+16q^8+24q^7+28q^6+27q^5+16q^4+4q^3)x^6+\\
(q^{16}+2q^{15}+5q^{14}+10q^{13}+18q^{12}+30q^{11}+47q^{10}+62q^9+76q^8+76q^7+62q^6+32q^5+8q^4)x^7 +O(x^8).
\end{multline*}
We refer to Figure~\ref{polyo14} for an illustration of the 14 polyominoes counted by the boldfaced coefficient of $x^4$.

The first terms of the series expansion of 
\begin{align*}F(x,u,q)=&\begin{pmatrix}
1 & 1 & 1 &  1
\end{pmatrix}\cdot \mathbf{F}(x,u,q)\\
=&F_{00}(x,u,q)+F_{01}(x,u,q)+F_{10}(x,u,q)+F_{11}(x,u,q)\\
=&G_{00}(x,u^2,q)+G_{01}(x,u^2,q)/u+G_{10}(x,u^2,q)+G_{11}(x,u^2,q)/u\end{align*}
are 
\begin{multline*}
qux+uq(qu+1)x^2+uq^2(q^2u^2+qu+q+u+1)x^3+\\
uq^2(q^4u^3+q^3u^2+q^3u+q^2u^2+2q^2u+qu^2+q^2+2qu+2q+2)x^4+\\
uq^3\left(q^6u^4+q^5u^3+q^5u^2+q^4u^3+2q^4u^2+q^3u^3+q^4u+3q^3u^2+\right.\\
\left.q^2u^3+q^4+2q^3u+3q^2u^2+2q^3+4q^2u+4q^2+5qu+5q+2u+2\right)x^5+O(x^6).
\end{multline*}

Finally the first terms of the series expansion of 
$F(1,1,q)=G(1,1,q)$ are 
$$2q+5q^2+15q^3+47q^4+149q^5+473q^6+1506q^7+4798q^8+O(q^9).$$
This sequence of coefficients does not appear in the On-Line Encyclopedia of Integer Sequences \cite{OEIS}.

Moreover, the series expansion of each $F_{ab}(x,u,q)$, for $a,b\in\{0,1\}$, can be extracted from the full series by retaining only those monomials $u^nx^m$ with $m\mod{2}=a$ and $n\mod{2}=b$. 

\section{A matrix continued fraction approach using automaton}

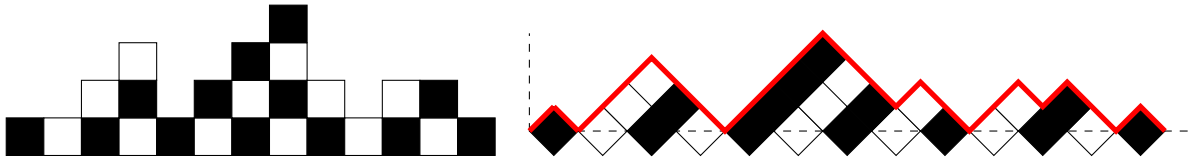
\begin{figure}[ht!]
\centering
\begin{tikzpicture}[scale=0.5]
\filldraw[fill=black] (0,0) rectangle ++(1,1);
\foreach \x in {1,2,3,4,5,6,7,8,9,10,11,12} {
    \case{\x}{0}
}
\foreach \y in {1} {
    \case{2}{\y}
}
\foreach \y in {1,2} {
    \case{3}{\y}
}
\foreach \y in {1} {
    \case{5}{\y}
}
\foreach \y in {1,2} {
    \case{6}{\y}
}
\foreach \y in {1,2,3} {
    \case{7}{\y}
}
\foreach \y in {1} {
    \case{8}{\y}
}
\foreach \y in {1} {
    \case{10}{\y}
}
\foreach \y in {1} {
    \case{11}{\y}
}
\end{tikzpicture}\quad
\begin{tikzpicture}[scale=0.65]
\newcommand{\lon}[2]{%
  \begin{scope}[shift={(#1,#2)}]
    \filldraw[fill=black, draw=black] (0,0) -- (0.5,0.5) -- (1,0) -- (0.5,-0.5) -- cycle;
  \end{scope}
  }
\newcommand{\lob}[2]{%
  \begin{scope}[shift={(#1,#2)}]
    \filldraw[fill=white, draw=black] (0,0) -- (0.5,0.5) -- (1,0) -- (0.5,-0.5) -- cycle;
  \end{scope}
}
 \lon{0}{0}
 \lon{2}{0}
 \foreach \y in {1,3,5,7,9,11} {
    \lob{\y}{0}
}
 \foreach \y in {2,4,6,8,10,12} {
    \lon{\y}{0}
}
 \foreach \y in {1.5,5.5,7.5,9.5} {
    \lob{\y}{0.5}
}
 \foreach \y in {2.5,4.5,6.5,10.5} {
    \lon{\y}{0.5}
}
\lon{5}{1}\lon{5.5}{1.5}
\draw[dashed] (0,0)--(13.5,0);
\draw[dashed] (0,0)--(0,2);
\draw[line width=0.7mm,red] (0,0)--(0.5,0.5);
\draw[line width=0.7mm,red] (0.5,0.5)--(1,0);
\draw[line width=0.7mm,red] (1,0)--(2.5,1.5)--(4,0)--(6,2)--(7.5,0.5)--(8,1)--(9,0)--(9.5,0.5)--(10,1)--(10.5,0.5)--(11,1)--(12,0)--(12.5,0.5)--(13,0);
\end{tikzpicture}
\caption{The black cell capacity of a Catalan polyomino is interpreted via Dyck paths.}
\label{CapacityDyck}
\end{figure}

In this section we work in the setting of Dyck paths. Recall that, from a Dyck path of semilength $n$, there exists a unique polyomino with 
$n$ columns such that the $i$th column contains $h$ cells, where 
$h$ is the height of the endpoint of the $i$th up-step of the path.
 In this context, it is well known that the area of a Dyck path  $P$ is the number of cells lying below the path and above the line $y=-1$, which also coincides with the area of the polyomino associated with $P$. One can also observe that the black cell capacity of the corresponding polyomino is the number of cells lying below the path and above the line $y=-1$ whose center $(a,b)$ satisfies $a-b=1 \mod{4}$. See Figure~\ref{CapacityDyck} for an illustration of the black capacity on a Dyck path.

We construct an automaton that generates partial Dyck paths. The states are organized into two layers, determined by the parity of the current height and the number of up-steps. A (red or black) arrow from state $i$ to $i+1$ represents the addition of an up-step at the end of the current path. The label on this arrow indicates the number of black cells created by this operation. An arrow from state $i+1$ to $i$ represents the addition of a down-step at the end of the current path. In this case, the operation does not create any black cells. For further examples of how infinite automata can be used to enumerate lattice paths, see, for instance, \cite{DCastro, Prod}.

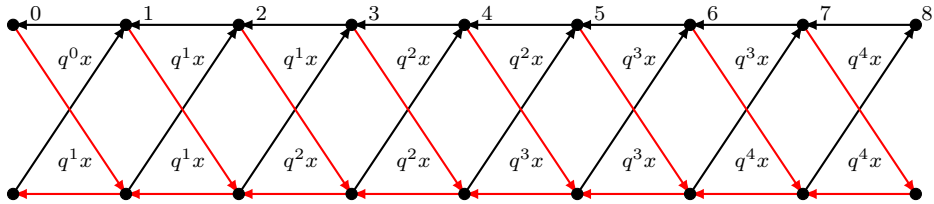
\begin{figure}[ht!]
\begin{center}
 			\begin{tikzpicture}[scale=1.5,main node/.style={circle,draw,font=\Large\bfseries}]

 				\foreach \x in {0,...,7}
 				{				
 					\draw[thick,-latex ] (\x,-1.5) to (\x+1,0) ;	
                    	\draw[thick, latex-] (\x,0) to  (\x+1,0);	
 							\draw[thick,red, -latex] (\x,0) to  (\x+1,-1.5);	
 				\draw[thick,red, latex-] (\x,-1.5) to  (\x+1,-1.5);	
                \node at  (\x+0.2,0.1){\tiny$\x$};
          \pgfmathtruncatemacro{\cx}{ceil((\x+1)/2)}
              \node at  (\x+0.55,-1.2){\tiny$q^{\cx}x$};
                        \pgfmathtruncatemacro{\fx}{floor((\x+1)/2)}
              \node at  (\x+0.55,-0.3){\tiny$q^{\fx}x$};
 	 				}

 				\node at  (8+0.1,0.1){\tiny$8$};
 				
 				\foreach \x in {0,1,2,3,4,5,6,7,8}
 				{
 					\draw (\x,0) circle (0.05cm);
 					\fill (\x,0) circle (0.05cm);
 					\draw (\x,-1.5) circle (0.05cm);
 					\fill (\x,-1.5) circle (0.05cm);
 			 				}
 			\end{tikzpicture}
 		\end{center}
 		\caption{Graph (automaton) to recognize  partial Dyck paths with respect to the length (marked with $x$) and the black cell capacity (marked with $q$). }
 		\label{automaton}
 	\end{figure}

From the automaton we obtain two families of generating functions,
$\{A_h(x,q)\}_{h\geq 0}$ and $\{B_h(x,q)\}_{h\geq 0}$, corresponding to the two
layers of states. For brevity, we write $A_h(x,q)=A_h$ and $B_h(x,q)=B_h$.
The recurrences follow by decomposing  the corresponding paths (equivalently, the associated polyominoes) according to their last step  (equivalently, their last column).

For each $h\geq 1$, we have
\begin{equation*}
A_h = A_{h-1} + a_h B_{h+1},
\qquad
B_h = B_{h-1} + b_zh A_{h+1},
\end{equation*}
where $a_h := xq^{\left\lceil\frac{h+1}{2}\right\rceil}$ and 
$b_h := x q^{\left\lfloor\frac{h+1}{2}\right\rfloor}$. Moreover, the initial conditions are
\begin{equation*}
A_0 = 1 + a_0 B_1,
\qquad
B_0 = 1 + b_0 A_1.
\end{equation*}

Our main goal is to obtain an explicit expression for $A_0$ (and $B_0$) in terms
of a matrix continued fraction (cf.\ Raissouli--Kacha \cite{RaissouliKacha2000} and Sorokin \cite{sorokin}).

Let
\[
\mathbf{U}_h:=\begin{pmatrix}0&a_h\\ b_h&0\end{pmatrix}
\qquad (h\geq 0),
\]
 and let $\mathbf{I}$ denote the $2\times2$
identity matrix. We consider the infinite matrix continued fraction
\begin{equation}\label{eq:S0-CF}
\mathbf{S}_0
=
\mathbf{I}-\cfrac{\mathbf{U}_0}{\,\mathbf{I}-\cfrac{\mathbf{U}_1}{\,\mathbf{I}-\cfrac{\mathbf{U}_2}{\,\mathbf{I}-\ddots}}}, \qquad\text{where }\ \cfrac{\mathbf{U}}{\mathbf{M}}:=\mathbf{U}\cdot\mathbf{M}^{-1}. 
\end{equation}
We define the $N$th convergent of \eqref{eq:S0-CF} by truncation:
\begin{equation}\label{eq:S0-conv}
\mathbf{S}^{(N)}_0
=
\mathbf{I}-\cfrac{\mathbf{U}_0}{\,\mathbf{I}-\cfrac{\mathbf{U}_1}{\,\cfrac{\ddots}{\mathbf{I}-\cfrac{\mathbf{U}_N}{\mathbf{I}}}}}\qquad (N\geq 0).
\end{equation}
By \eqref{eq:S0-CF} we mean the unique $2\times 2$ matrix $\mathbf{S}_0$ such that for every $k\geq 1$ there exists $N_k$ with
\[
\mathbf{S}^{(N)}_0 \equiv \mathbf{S}_0 \pmod{x^k}\qquad\text{for all }N\geq N_k.
\]
Equivalently, for each  entry, the formal power series expansions in $x$ of $\mathbf{S}^{(N)}_0$ and
$\mathbf{S}_0$ agree in every term of degree $<k$ (that is, they coincide up to order $x^{k-1}$).

For example, expanding the first few  convergents $\mathbf{S}_0$ we obtain
\begin{align*}
\mathbf{S}_0&=\mathbf{I}-\cfrac{\begin{pmatrix} 0 & q x \\
 x & 0 \\\end{pmatrix}}{\,\mathbf{I}-\cfrac{\begin{pmatrix}  0 & q x \\
 q x & 0 \\\end{pmatrix}}{\,\mathbf{I}-\cfrac{\begin{pmatrix}  0 & q^2 x \\
 q x & 0 \\\end{pmatrix}}{\,\mathbf{I}-\ddots}}} \\
&=
\mathbf{I}
+ x\,\mathbf{M}_1
+ x^2\,\mathbf{M}_2
+ x^3\,\mathbf{M}_3
+ x^4\,\mathbf{M}_4
+ x^5\,\mathbf{M}_5
+ O(x^6),
\end{align*}
where 
\begin{align*}
\mathbf{M}_1&=
\begin{pmatrix}
0 & -q\\
-1& 0
\end{pmatrix},
\qquad
\mathbf{M}_2=
\begin{pmatrix}
-q^2 & 0\\
0 & -q
\end{pmatrix},
\qquad
\mathbf{M}_3=
\begin{pmatrix}
0 & -(q^4+q^3)\\
-2q^2 & 0
\end{pmatrix},\\
\mathbf{M}_4&=
\begin{pmatrix}
-(q^6+2q^5+2q^4) & 0\\
0 & -(3q^4+2q^3)
\end{pmatrix},\\
\mathbf{M}_5&=
\begin{pmatrix}
0 & -(q^9+2q^8+4q^7+5q^6+2q^5)\\
-(4q^6+6q^5+4q^4) & 0
\end{pmatrix}.
\end{align*}

\begin{lemma}\label{lem:existenceS0}
For each $N\geq 0$ the convergent $\mathbf{S}_0^{(N)}$ in \eqref{eq:S0-conv} is well-defined, that is, every matrix occurring in a denominator is invertible; in particular, $\mathbf{S}_0^{(N)}$ is invertible. Moreover, the entries of $\mathbf{S}_0^{(N)}$ stabilize coefficientwise in $x$, and therefore
\eqref{eq:S0-CF} defines a unique $2\times 2$ matrix $\mathbf{S}_0$ such that
\[
\mathbf{S}_0^{(N)} \equiv \mathbf{S}_0 \pmod{x^k}
\quad\text{for $N$ large enough.}
\]
\end{lemma}

\begin{proof}
Write $\mathcal{R}=\Q(q)[[x]]$. Since $a_h,b_h\in x\cdot\mathcal{R}$, we have $\mathbf{U}_h\in x\cdot\mathcal{M}_2(\mathcal{R})$ for all $h\geq 0$, where $\mathcal{M}_2(\mathcal{R})$ is the set of $2\times 2$ matrices on $\mathcal{R}$. 
Define the truncations recursively by setting $\mathbf{S}_{N+1}^{(N)}:=\mathbf{I}$ and
\[
\mathbf{S}_h^{(N)}:=\mathbf{I}-\mathbf{U}_h\bigl(\mathbf{S}_{h+1}^{(N)}\bigr)^{-1}\qquad(0 \leq h \leq N),
\]
so that $\mathbf{S}_0^{(N)}$ coincides with \eqref{eq:S0-conv}.
Since $\mathbf{S}_{N+1}^{(N)}=\mathbf{I}\equiv \mathbf{I}\pmod{x}$, a backward induction gives
$\mathbf{S}_h^{(N)}\equiv \mathbf{I}\pmod{x}$ for all $0\leq h\leq N$. Hence each $\mathbf{S}_h^{(N)}$
is invertible in $\mathcal{M}_2(\mathcal{R})$.

It remains to justify that the coefficients of $\mathbf{S}_0^{(N)}$ stabilize as $N\to\infty$.
Fix $k\geq 1$. Since $\mathbf{U}_h\in x\cdot\mathcal{M}_2(\mathcal{R})$, every additional level of the continued fraction
introduces at least one extra factor of $x$. Consequently, the truncation error at depth $N$
starts in degree $x^{N+1}$, and in particular $\mathbf{S}_0^{(N+1)}\equiv \mathbf{S}_0^{(N)}\pmod{x^k}$ for all $N\geq k-1$.

Therefore the entries of $\mathbf{S}_0^{(N)}$ stabilize coefficientwise in $x$: for each fixed $j\geq 0$, the coefficient of $x^j$ in $\mathbf{S}_0^{(N)}$ eventually becomes independent of $N$. We then define $\mathbf{S}_0\in \mathcal{M}_2(\mathcal{R})$ by taking, for each $j$, the coefficient of $x^j$ in $\mathbf{S}_0$ to be this eventual value. In particular, since $\mathbf{S}_0^{(N)}\equiv \mathbf{I}\pmod{x}$ for all $N$, the constant term of $\mathbf{S}_0$ is $\mathbf{I}$, and hence $\mathbf{S}_0\equiv \mathbf{I}\pmod{x}$.
\end{proof}

\begin{theorem}\label{thm:A0B0}
Let $\mathbf{S}_0$ be the matrix defined by the infinite continued fraction \eqref{eq:S0-CF},
and let $\mathbf{e}=\begin{pmatrix} 1\\1\end{pmatrix}$. Then 
\[
\begin{pmatrix}
    A_0\\ B_0
\end{pmatrix}=\mathbf{S}_0^{-1}\mathbf{e},
\]
that is,
\[
A_0=\begin{pmatrix}1&0\end{pmatrix}\mathbf{S}_0^{-1}\begin{pmatrix}
    1\\ 1
\end{pmatrix},
\qquad
B_0=\begin{pmatrix}0&1\end{pmatrix}\mathbf{S}_0^{-1}\begin{pmatrix}
    1\\ 1
\end{pmatrix}.
\]
\end{theorem}
\begin{proof}
Set
\[
\mathbf{V}_h:=\begin{pmatrix}A_h\\ B_h\end{pmatrix},
\qquad
\mathbf{e}:=\begin{pmatrix}1\\1\end{pmatrix},
\qquad
\mathbf{U}_h:=\begin{pmatrix}0&a_h\\ b_h&0\end{pmatrix}.
\]
Then the given system is equivalently
\begin{equation}\label{eq:Vsystem}
\mathbf{V}_0=\mathbf{e}+\mathbf{U}_0\mathbf{V}_1,
\qquad
\mathbf{V}_h=\mathbf{V}_{h-1}+\mathbf{U}_h\mathbf{V}_{h+1}\quad(h\geq 1).
\end{equation}
For each $h\geq 0$, let $\mathbf{S}_h$ denote the matrix continued fraction
\[
\mathbf{S}_h
=
\mathbf{I}-\cfrac{\mathbf{U}_h}{\,\mathbf{I}-\cfrac{\mathbf{U}_{h+1}}{\,\mathbf{I}-\cfrac{\mathbf{U}_{h+2}}{\,\mathbf{I}-\ddots}}}\,,
\]
so that $\mathbf{S}_0$ is the matrix in \eqref{eq:S0-CF}.  Notice that by Lemma~\ref{lem:existenceS0}, each
$\mathbf{S}_h$ is well-defined and invertible, moreover $\mathbf{S}_h\equiv \mathbf{I}\pmod{x}$.

Define the sequence $(W_h)_{h\geq 0}$ by
\[
\mathbf{W}_0:=\mathbf{S}_0^{-1}\mathbf{e},
\qquad
\mathbf{W}_{h+1}:=\mathbf{S}_{h+1}^{-1}\mathbf{W}_h \quad (h\geq 0).
\]
Equivalently, $\mathbf{W}_{h-1}=\mathbf{S}_h\mathbf{W}_h$ for all $h\geq 1$. 

We claim that $(\mathbf{W}_h)_{h\geq 0}$ satisfies the same system \eqref{eq:Vsystem}.  First, since
$\mathbf{S}_0=\mathbf{I}-\mathbf{U}_0\mathbf{S}_1^{-1}$ (this is just the first step of the continued fraction), we get
\[
\mathbf{e}=\mathbf{S}_0\mathbf{W}_0=(\mathbf{I}-\mathbf{U}_0\mathbf{S}_1^{-1})\mathbf{W}_0=\mathbf{W}_0-\mathbf{U}_0\mathbf{W}_1,
\]
hence $\mathbf{W}_0=\mathbf{e}+\mathbf{U}_0\mathbf{W}_1$, which is the $h=0$ equation in \eqref{eq:Vsystem}.
Next, for $h\geq 1$, using $\mathbf{S}_h=\mathbf{I}-\mathbf{U}_h\mathbf{S}_{h+1}^{-1}$ and $\mathbf{W}_{h+1}=\mathbf{S}_{h+1}^{-1}\mathbf{W}_h$ we obtain
\[
\mathbf{W}_{h-1}=\mathbf{S}_h\mathbf{W}_h=(\mathbf{I}-\mathbf{U}_h\mathbf{S}_{h+1}^{-1})\mathbf{W}_h=\mathbf{W}_h-\mathbf{U}_h\mathbf{W}_{h+1},
\]
so,
\[
\mathbf{W}_h=\mathbf{W}_{h-1}+\mathbf{U}_h\mathbf{W}_{h+1}\qquad(h\geq 1),
\]
which is the second equation in \eqref{eq:Vsystem}.

Therefore $(\mathbf{W}_h)$ is a solution of the original system. By uniqueness of the solution, we must have $\mathbf{W}_h=\mathbf{V}_h$ for all $h$, and in particular
\[
\begin{pmatrix}A_0\\ B_0\end{pmatrix}=\mathbf{V}_0=\mathbf{W}_0=\mathbf{S}_0^{-1}\mathbf{e}.
\]
\end{proof}

Recall that $\mathbf{S}_0$ satisfies $\mathbf{S}_0\equiv \mathbf{I}\pmod{x}$, hence we may write
\[
\mathbf{S}_0=\mathbf{I}+\mathbf{X} \qquad\text{with } \mathbf{X}\in x\cdot\mathcal{M}_2(\Q(q)[[x]]).
\]
In particular, $\mathbf{X}$ has no constant term in $x$, so the inverse of $\mathbf{S}_0$ can be computed
formally by the geometric-series identity
\[
\mathbf{S}_0^{-1}=(\mathbf{I}+\mathbf{X})^{-1}=\mathbf{I}-\mathbf{X}+\mathbf{X}^2-\mathbf{X}^3+\cdots,
\]
which is well-defined coefficientwise in $x$ because $\mathbf{X}^m$ starts in degree $x^m$.

Using the previous expansion of $\mathbf{S}_0$, we obtain
\[
\mathbf{S}_0^{-1}
=
\mathbf{I}
+ x\,\mathbf{N}_1
+ x^2\,\mathbf{N}_2
+ x^3\,\mathbf{N}_3
+ x^4\,\mathbf{N}_4
+ x^5\,\mathbf{N}_5
+ O(x^6),
\]
where
\begin{align*}
\mathbf{N}_1&=
\begin{pmatrix}
0 & q\\
1& 0
\end{pmatrix},
\qquad
\mathbf{N}_2=
\begin{pmatrix}
q(q+1) & 0\\
0 & 2q
\end{pmatrix},
\qquad
\mathbf{N}_3=
\begin{pmatrix}
0 & q^2(q^2+2q+2)\\
q(3q+2)& 0
\end{pmatrix},\\
\mathbf{N}_4&=
\begin{pmatrix}
q^2(q^4+2q^3+4q^2+5q+2) & 0\\
0 & q^2(4q^2+6q+4)
\end{pmatrix},\\
\mathbf{N}_5&=
\begin{pmatrix}
0 & q^3(q^6+2q^5+5q^4+8q^3+12q^2+10q+4)\\
q^2(5q^4+8q^3+13q^2+12q+4) & 0
\end{pmatrix}.
\end{align*}

Notice that 
$$\begin{pmatrix}1&0\end{pmatrix}\mathbf{N}_4\begin{pmatrix}
    1\\ 1  
\end{pmatrix}=q^6+2q^5+4q^4+5q^3+2q^2,$$
which corresponds to the distribution of polyominoes of length four with respect to the black cell capacity (see Figure~\ref{polyo14} for an illustration of these polyominoes).

\begin{remark} The methods employed in these last two sections enable us to derive the series expansions of the above generating functions. However, they do not provide closed-form expressions for these functions. In the next section, we exploit the relationship between the black capacity and the vertical black capacity established in Theorem 1 in order to obtain a closed-form expression.
\end{remark}

\section{Length, vertical black/white cell capacity and last value}
In this section we count Catalan polyominoes by their length (number of columns), by their
vertical black/white cell capacities, and by the height of the last column. Let us recall that the vertical black (resp. vertical white) cell capacity is the total number of cells in the columns of odd indices (resp. even indices). 

Let $EV(x,y,z,u)$ (resp. $OD(x,y,z,u)$) be the generating function where the coefficient $ev_{n,i,j,k}$ (resp. $od_{n,i,j,k}$) of $x^ny^iz^ju^k$ is the number of Catalan polyominoes $P$ with $n$ columns,  $n$ even (resp. odd), such that  $\ver(P)=i$, $\white(P)=j$, and $\last(P)=k$.

\begin{theorem} The generating functions $EV(x,y,z,u)$ and $OD(x,y,z,u)$ satisfy the following system of equations:
   \[\begin{cases} 
   \begin{aligned}EV(x,y,z,u)&=\frac{xzu}{zu-1}\left(uzOD(x,y,z,zu)-OD(x,y,z,1)\right),\\
    OD(x,y,z,u)&=xyu+\frac{xyu}{yu-1}\left(uyEV(x,y,z,yu)-EV(x,y,z,1)\right),
\end{aligned}
\end{cases}\]
which is equivalent to the following matrix equation:
 \begin{align}
 \mathbf{V}(u)=\mathbf{M}(u)\cdot \mathbf{V}(yzu)-\mathbf{N}(u)\cdot \mathbf{V}(1)+\mathbf{B}(u),
 \end{align} 
 where
\[\mathbf{M}(u)=
\begin{pmatrix}
\frac{x^2y^2z^4u^4}{(uz-1)(uyz-1)}& 0 \\ 
 0& \frac{x^2y^4z^2u^4}{(uy-1)(uyz-1)}
\end{pmatrix},
 \quad \mathbf{N}(u)=
\begin{pmatrix}
\frac{x^2yz^3u^3}{(uz-1)(uyz-1)}&  \frac{xzu}{uz-1}\\
\frac{xyu}{uy-1}& \frac{x^2y^3zu^3}{(uy-1)(uyz-1)}
\end{pmatrix},
\]

\[
\mathbf{B}(u)=
\begin{pmatrix} \frac{x^2yz^3u^3}{uz-1}\\
                xyu
\end{pmatrix}\quad \mbox{ and } \quad \mathbf{V}(u)=
\begin{pmatrix} EV(x,y,z,u)\\
                OD(x,y,z,u)
\end{pmatrix}.
\]
\end{theorem}
\begin{proof}
We distinguish two types of Catalan polyominoes $P$ according to the parities of the number  of columns \lth(P).

Case 1: $\lth(P)$ is odd. If  $P$ has a single column, it consists of one cell, and its contribution to the generating function is $xyu$. Otherwise, assume that the polyomino has at least two columns. $P$ is obtained from a polyomino $Q$ of length $\lth(Q)=\lth(P)-1=0\mod{2}$ by adding on its right a column with $\last(P)$ cells, by preserving the Catalan structure, i.e. the condition $\last(Q)+1\geq \last(P)$. 
    Then the contribution for these polyominoes is 
\begin{align*}
    A_1:=&x\sum_{n,i,j,k} ev_{2n,i,j,k}\,x^{2n}y^iz^j\left(uy+ u^2y^2+\cdots + u^{k+1}y^{k+1}\right).\\
    =&xuy\sum_{n,i,j,k} ev_{2n,i,j,k}\,x^{2n}y^iz^j\left(\frac{u^{k+1}y^{k+1}-1}{uy-1}\right)\\
    =&\frac{xuy}{uy-1}\left(uyEV(x,y,z,u)-EV(x,y,z,1)\right), 
\end{align*}
which gives the second equation.
 The second case ($\lth(P)$ even) gives the first equation  {\it mutatis mutandis}.
\end{proof}

By iterating Equation (5), we obtain a formal expression for the vector $\mathbf{V}(u)$.

\begin{theorem} We have
 \[
\mathbf{V}(u) = \sum_{k=0}^{\infty} \mathbf{P}_k(u) \cdot \Big( \mathbf{B}((yz)^k u) - \mathbf{N}((yz)^k u) \, \mathbf{V}(1) \Big),
\]
with 
\[
\mathbf{P}_k(u) := \prod_{j=0}^{k-1} \mathbf{M}((yz)^j u), \quad \mbox{ and } \quad \mathbf{P}_{-1}(u)=\mathbf{P}_0(u)= \mathbf{I},
\]
where $\mathbf{I}$ is the identity matrix of size 2.
\end{theorem}

Substituting $u=1$ in the expression of $\mathbf{V}(u)$, we obtain:
\[
\mathbf{V}(1) = \sum_{k=0}^{\infty} \mathbf{P}_k(1) \cdot\Big( \mathbf{B}((yz)^k) - \mathbf{N}((yz)^k) \, \mathbf{V}(1) \Big).
\]
Isolating $\mathbf{V}(1)$ we obtain
\[
\left( \mathbf{I} + \sum_{k=0}^{\infty} \mathbf{P}_k(1) \mathbf{N}((yz)^k) \right) \cdot\mathbf{V}(1) = \sum_{k=0}^{\infty} \mathbf{P}_k(1) \mathbf{B}((yz)^k),
\]
and we deduce the following theorem.

\begin{theorem} We have 
\[
 \mathbf{V}(1)=\left( \mathbf{I} + \mathbf{S}(x,y,z) \right)^{-1}\cdot\sum_{k=0}^{\infty} \mathbf{P}_k(1) \mathbf{B}((yz)^k),
 \]
with
\[\mathbf{S}(x,y,z)=\sum_{k=0}^{\infty} \mathbf{P}_k(1) \mathbf{N}((yz)^k) \quad\mbox{ and }\quad
\mathbf{P}_k(1)=
\begin{pmatrix}
G_k(y,z)
& 0\\
0 &G_k(z,y)
\end{pmatrix}
,\] 
with \[G_k(y,z)=\frac{x^{2k} y^{2k} z^{4k}  (yz)^{2k(k-1)}}
{( z;\,yz)_k\,( y z;\,yz)_k},\] 
 and $(a;\,b)_k$ is the Pochhammer symbol (see~\cite{gasper}) 
 \[
(a;\,b)_k
=
\prod_{j=0}^{k-1} \bigl(1-a b^j\bigr). 
\]
\end{theorem}

Now let us show how we can calculate $\mathbf{V}(1)$. If we set 
\[ A_k(y,z)=\frac{x^2yz^3(yz)^{3k}}{((yz)^kz-1)((yz)^kyz-1)} \quad \mbox{ and }\quad B_k(y,z)=\frac{xz(yz)^k}{(yz)^kz-1},
\]
then we have \[\mathbf{P}_k(1)\mathbf{N}((yz)^k)=
\begin{pmatrix}
G_k(y,z)A_k(y,z)
& G_k(y,z)B_k(y,z)\\
G_k(z,y)B_k(z,y) &G_k(z,y)A_k(z,y)
\end{pmatrix},\]
which implies that
\[
\mathbf{S}(x,y,z)=
\begin{pmatrix}
\phi(y,z)
& \psi(y,z)\\
\psi(z,y) &\phi(z,y)
\end{pmatrix},
\]
with 
\[\phi(y,z)=
\sum\limits_{k=0}^\infty\frac{x^{2k+2} y^{2k+1} z^{4k+3} (yz)^{k(2k+1)}}
{(z;\,yz)_{k+1}\,(y z;\,yz)_{k+1}} \mbox{ and }
\psi(y,z)=
-\sum\limits_{k=0}^\infty\frac{x^{2k+1} y^{2k} z^{4k+1} (yz)^{k(2k-1)}}
{(z;\,yz)_{k+1}\,(y z;\,yz)_k}.
\]
Now we compute the inverse of $\mathbf{I}+\mathbf{S}(x,y,z)$, 
\[
(\mathbf{I}+\mathbf{S}(x,y,z))^{-1}=
\frac{1}{\Delta}
\begin{pmatrix}
1+\phi(z,y)
& -\psi(y,z)\\
-\psi(z,y) &1+\phi(y,z)
\end{pmatrix} \]
where \[\Delta=(1+\phi(y,z))(1+\phi(z,y))-\psi(y,z)\psi(z,y)
.\]

Finally, we obtain a close form for $\mathbf{V}(1)$.
\begin{theorem} We have 
    \[
\mathbf{V}(1)=\frac{1}{\Delta}
\begin{pmatrix}
-(1+\phi(z,y))f(y,z)-\psi(y,z)g(y,z)\\
(1+\phi(y,z))g(y,z)+\psi(z,y)f(y,z)
\end{pmatrix} 
\]
with 

\[ 
f(y,z)=\sum_{k=0}^\infty
\frac{x^{2k+2} y^{2k+1} z^{4k+3} (yz)^{k(2k+1)}}
{(z;\,yz)_{k+1}\,(y z;\,yz)_k}, \quad
g(y,z)=\sum_{k=0}^\infty
\frac{x^{2k+1} y^{4k+1} z^{2k} (yz)^{k(2k-1)}}
{(y;\,yz)_k\,(y z;\,yz)_k},
\]
\[\phi(y,z)=
\sum\limits_{k=0}^\infty\frac{x^{2k+2} y^{2k+1} z^{4k+3} (yz)^{k(2k+1)}}
{(z;\,yz)_{k+1}\,(y z;\,yz)_{k+1}},\quad 
\psi(y,z)=
-\sum\limits_{k=0}^\infty\frac{x^{2k+1} y^{2k} z^{4k+1} (yz)^{k(2k-1)}}
{(z;\,yz)_{k+1}\,(y z;\,yz)_k},
\]
and 
\[\Delta=(1+\phi(y,z))(1+\phi(z,y))-\psi(y,z)\psi(z,y)
.\]
\end{theorem}
The first terms of the series expansion of the generating function  $OD(x,y,1,1)$ (the second coordinate of $\mathbf{V}(1)$ evaluated at $z=1$)
are  
\begin{align*}
&xy+(y^2+2y+2)y^2x^3+y^3(y^6+2y^5+5y^4+8y^3+12y^2+10y+4)x^5+\\
&(y^{12}+2y^{11}+5y^{10}+10y^9+18y^8+30y^7+47y^6+62y^5+76y^4+76y^3+62y^2+32y+8)y^4x^7 +O(x^9).
\end{align*}
Due to Theorem~\ref{th11}, these terms  correspond exactly to the odd terms of $F(x,1,y)$ (i.e. the terms in $x^{2n-1}$, $n\geq 1$) obtained in the previous section.

The first terms of the series expansion of the generating function  $EV(x,1,z,1)$ (the first coordinate of $\mathbf{V}(1)$ evaluated at $y=1$)
are  
\begin{align*}
&(z+1)zx^2+\bm{z^2(z^4+2z^3+4z^2+5z+2)x^4}+\\
&(z^9+2z^8+5z^7+9z^6+16z^5+24z^4+28z^3+27z^2+16z+4)z^3x^6+O(x^8).
\end{align*}
Due to Theorem~\ref{th11}, these terms  correspond exactly to the even terms of $F(x,1,z)$ (i.e. the terms in $x^{2n}$, $n\geq 1$) obtained in the previous section. We refer to Figure~\ref{polyo14} for an illustration of of the polyominoes of length 4  counted by the boldfaced coefficient of $x^4$.

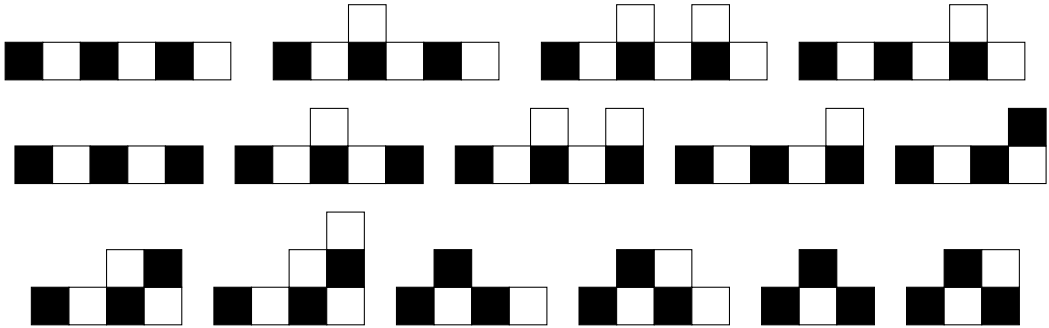
\begin{figure}[H]
\centering
\begin{tikzpicture}[scale=0.5]-
\filldraw[fill=black] (0,0) rectangle ++(1,1);
\foreach \x in {1,2,3,4,5} {
    \case{\x}{0}
}
\end{tikzpicture}
\quad\begin{tikzpicture}[scale=0.5]-
\filldraw[fill=black] (0,0) rectangle ++(1,1);
\foreach \x in {1,2,3,4,5} {
    \case{\x}{0}
}
\case{2}{1}
\end{tikzpicture}
\quad\begin{tikzpicture}[scale=0.5]-
\filldraw[fill=black] (0,0) rectangle ++(1,1);
\foreach \x in {1,2,3,4,5} {
    \case{\x}{0}
}
\case{2}{1}\case{4}{1}
\end{tikzpicture}\quad 
\begin{tikzpicture}[scale=0.5]-
\filldraw[fill=black] (0,0) rectangle ++(1,1);
\foreach \x in {1,2,3,4,5} {
    \case{\x}{0}
}
\case{4}{1}
\end{tikzpicture}\\[0.8em]
\quad \begin{tikzpicture}[scale=0.5]-
\filldraw[fill=black] (0,0) rectangle ++(1,1);
\foreach \x in {1,2,3,4} {
    \case{\x}{0}
}
\end{tikzpicture}\quad 
\begin{tikzpicture}[scale=0.5]-
\filldraw[fill=black] (0,0) rectangle ++(1,1);
\foreach \x in {1,2,3,4} {
    \case{\x}{0}
}
    \case{2}{1}
\end{tikzpicture}\quad 
\begin{tikzpicture}[scale=0.5]-
\filldraw[fill=black] (0,0) rectangle ++(1,1);
\foreach \x in {1,2,3,4} {
    \case{\x}{0}
}
    \case{2}{1}
    \case{4}{1}
\end{tikzpicture}\quad 
\begin{tikzpicture}[scale=0.5]-
\filldraw[fill=black] (0,0) rectangle ++(1,1);
\foreach \x in {1,2,3,4} {
    \case{\x}{0}
}
    \case{4}{1}
\end{tikzpicture}\quad 
\begin{tikzpicture}[scale=0.5]-
\filldraw[fill=black] (0,0) rectangle ++(1,1);
\foreach \x in {1,2,3} {
    \case{\x}{0}
}
    \case{3}{1}
\end{tikzpicture}\\[0.8em]\quad 
\begin{tikzpicture}[scale=0.5]-
\filldraw[fill=black] (0,0) rectangle ++(1,1);
\foreach \x in {1,2,3} {
    \case{\x}{0}
}
    \case{2}{1}
    \case{3}{1}
\end{tikzpicture}\quad 
\begin{tikzpicture}[scale=0.5]-
\filldraw[fill=black] (0,0) rectangle ++(1,1);
\foreach \x in {1,2,3} {
    \case{\x}{0}
}
    \case{2}{1}
\foreach \y in {1,2} {
    \case{3}{\y}
}
\end{tikzpicture}\quad
\begin{tikzpicture}[scale=0.5]-
\filldraw[fill=black] (0,0) rectangle ++(1,1);
\foreach \x in {1,2,3} {
    \case{\x}{0}
}
    \case{1}{1}
\end{tikzpicture}\quad 
\begin{tikzpicture}[scale=0.5]-
\filldraw[fill=black] (0,0) rectangle ++(1,1);
\foreach \x in {1,2,3} {
    \case{\x}{0}
}
    \case{1}{1}
    \case{2}{1}
\end{tikzpicture}\quad 
\begin{tikzpicture}[scale=0.5]-
\filldraw[fill=black] (0,0) rectangle ++(1,1);
\foreach \x in {1,2} {
    \case{\x}{0}
}
    \case{1}{1}
\end{tikzpicture}\quad 
\begin{tikzpicture}[scale=0.5]-
\filldraw[fill=black] (0,0) rectangle ++(1,1);
\foreach \x in {1,2} {
    \case{\x}{0}
}
    \case{1}{1}
    \case{2}{1}
\end{tikzpicture}
\caption{The 15 Catalan polyominoes $P$ with $\bck(P)=3$.} 
\label{polyo15}
\end{figure}

Now, if we calculate the series expansions of $OD(1,z,1,1)+EV(1,1,z,1)$
we obtain 
$$2z+5z^2+15z^3+47z^4+149z^5+473z^6+1506z^7+4798z^8+O(z^9).$$
Due to Theorem~\ref{th11}, these terms correspond to the series expansion of $F(1,1,q)$ obtained in the previous section. See Figure~\ref{polyo15} for an illustration of the 15 Catalan polyominoes having a black cell capacity equal to 3.

As a byproduct, we obtain the generating function for the polyominoes classified by the length and the vertical black cell capacity, i.e.  $OD(x,y,1,1)+EV(x,y,1,1)$, where the first terms of the series expansion are
\begin{multline*}yx + 2yx^2 + y^2(y^2 + 2y + 2)x^3  + (4y^4 + 6y^3 + 4y^2)x^4\\
 + (y^6 + 2y^5 + 5y^4 + 8y^3 + 12y^2 + 10y + 4)y^3x^5 +\\  (6y^9 + 10y^8 + 22y^7 + 28y^6 + 34y^5 + 24y^4 + 8y^3)x^6+O(x^7).
\end{multline*}

Finally, the generating function for the Catalan polyominoes enumerated by the vertical black cell capacity is $OD(1,y,1,1)+EV(1,y,1,1)$, and the first terms of the series expansion are
\[3y+6y^2+20y^3+63y^4+166y^5+342y^6+553y^7+734y^8 + O(y^9).\]
This sequence of coefficients does not appear in the On-Line Encyclopedia of Integer Sequences \cite{OEIS}.

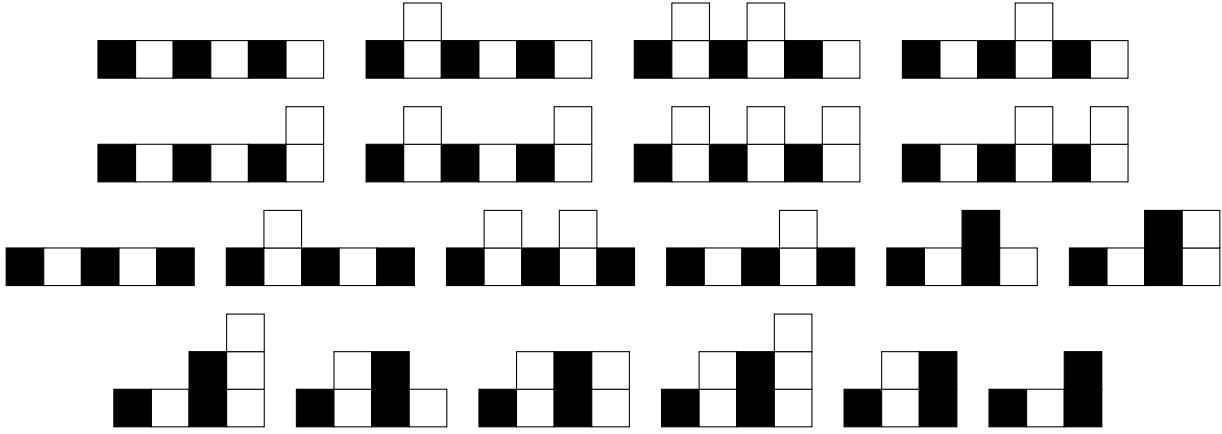
\begin{figure}[H]
\centering
\begin{tikzpicture}[scale=0.5]-
\filldraw[fill=black] (0,0) rectangle ++(1,1);
\foreach \x in {1,2,3,4,5} {
    \case{\x}{0}
}
\end{tikzpicture}
\quad\begin{tikzpicture}[scale=0.5]-
\filldraw[fill=black] (0,0) rectangle ++(1,1);
\foreach \x in {1,2,3,4,5} {
    \case{\x}{0}
}
\casee{1}{1}
\end{tikzpicture}
\quad\begin{tikzpicture}[scale=0.5]-
\filldraw[fill=black] (0,0) rectangle ++(1,1);
\foreach \x in {1,2,3,4,5} {
    \case{\x}{0}
}
\casee{1}{1}\casee{3}{1}
\end{tikzpicture}
\quad\begin{tikzpicture}[scale=0.5]-
\filldraw[fill=black] (0,0) rectangle ++(1,1);
\foreach \x in {1,2,3,4,5} {
    \case{\x}{0}
}
\casee{3}{1}
\end{tikzpicture}\\[0.8em]
\begin{tikzpicture}[scale=0.5]-
\filldraw[fill=black] (0,0) rectangle ++(1,1);
\foreach \x in {1,2,3,4,5} {
    \case{\x}{0}
}\casee{5}{1}
\end{tikzpicture}
\quad\begin{tikzpicture}[scale=0.5]-
\filldraw[fill=black] (0,0) rectangle ++(1,1);
\foreach \x in {1,2,3,4,5} {
    \case{\x}{0}
}
\casee{1}{1}\casee{5}{1}
\end{tikzpicture}
\quad\begin{tikzpicture}[scale=0.5]-
\filldraw[fill=black] (0,0) rectangle ++(1,1);
\foreach \x in {1,2,3,4,5} {
    \case{\x}{0}
}
\casee{1}{1}\casee{3}{1}\casee{5}{1}
\end{tikzpicture}
\quad\begin{tikzpicture}[scale=0.5]-
\filldraw[fill=black] (0,0) rectangle ++(1,1);
\foreach \x in {1,2,3,4,5} {
    \case{\x}{0}
}
\casee{3}{1}\casee{5}{1}
\end{tikzpicture}\\[0.8em]
\begin{tikzpicture}[scale=0.5]-
\filldraw[fill=black] (0,0) rectangle ++(1,1);
\foreach \x in {1,2,3,4} {
    \case{\x}{0}
}
\end{tikzpicture}\quad 
\begin{tikzpicture}[scale=0.5]-
\filldraw[fill=black] (0,0) rectangle ++(1,1);
\foreach \x in {1,2,3,4} {
    \case{\x}{0}
}
    \casee{1}{1}
\end{tikzpicture}\quad 
\begin{tikzpicture}[scale=0.5]-
\filldraw[fill=black] (0,0) rectangle ++(1,1);
\foreach \x in {1,2,3,4} {
    \case{\x}{0}
}
    \casee{1}{1}
    \casee{3}{1}
\end{tikzpicture}\quad 
\begin{tikzpicture}[scale=0.5]-
\filldraw[fill=black] (0,0) rectangle ++(1,1);
\foreach \x in {1,2,3,4} {
    \case{\x}{0}
}
    \casee{3}{1}
\end{tikzpicture}\quad 
\begin{tikzpicture}[scale=0.5]-
\filldraw[fill=black] (0,0) rectangle ++(1,1);
\foreach \x in {1,2,3} {
    \case{\x}{0}
}
    \casee{2}{1}
\end{tikzpicture}\quad 
\begin{tikzpicture}[scale=0.5]-
\filldraw[fill=black] (0,0) rectangle ++(1,1);
\foreach \x in {1,2,3} {
    \case{\x}{0}
}
    \casee{2}{1}
    \casee{3}{1}
\end{tikzpicture}\\[0.8em]
\begin{tikzpicture}[scale=0.5]-
\filldraw[fill=black] (0,0) rectangle ++(1,1);
\foreach \x in {1,2,3} {
    \case{\x}{0}
}
    \casee{2}{1}
    \casee{3}{1}
    \case{3}{2}
\end{tikzpicture}\quad
\begin{tikzpicture}[scale=0.5]-
\filldraw[fill=black] (0,0) rectangle ++(1,1);
\foreach \x in {1,2,3} {
    \case{\x}{0}
}
    \casee{2}{1}\casee{1}{1}
\end{tikzpicture}\quad 
\begin{tikzpicture}[scale=0.5]-
\filldraw[fill=black] (0,0) rectangle ++(1,1);
\foreach \x in {1,2,3} {
    \case{\x}{0}
}\casee{1}{1}
    \casee{2}{1}
    \casee{3}{1}
\end{tikzpicture}\quad 
\begin{tikzpicture}[scale=0.5]-
\filldraw[fill=black] (0,0) rectangle ++(1,1);
\foreach \x in {1,2,3} {
    \case{\x}{0}
}\casee{1}{1}
    \casee{2}{1}
    \casee{3}{1}
    \case{3}{2}
\end{tikzpicture}\quad\begin{tikzpicture}[scale=0.5]-
\filldraw[fill=black] (0,0) rectangle ++(1,1);
\foreach \x in {1,2} {
    \case{\x}{0}
}
    \casee{1}{1}
    \casee{2}{1}
\end{tikzpicture}\quad 
\begin{tikzpicture}[scale=0.5]-
\filldraw[fill=black] (0,0) rectangle ++(1,1);
\foreach \x in {1,2} {
    \case{\x}{0}
}
    \casee{2}{1}
\end{tikzpicture}
\caption{The 20 Catalan polyominoes $P$ with $\ver(P)=3$.} 
\label{polyovert}
\end{figure}

We refer to Figure~\ref{polyovert} for an illustration of the 20 polyominoes with $\ver(P)=3$.

\section{Other functional equation for the black cell capacity} 
We introduce the statistic $\texttt{s}$ on the Catalan polyominoes by setting $\texttt{s}(P)=\ver(P)$ when $\lth(P)$ is odd, and $\texttt{s}(P)=\white(P)$, otherwise. By Theorem~\ref{th11}, the statistics $\texttt{s}$ and $\bck$ are equidistributed on Catalan polyominoes.  We also define the complementary statistic $\bar{\texttt{s}}$ by  $\bar{\texttt{s}}(P)=\ver(P)$ when $\lth(P)$ is even, and $\bar{\texttt{s}}(P)=\white(P)$, otherwise.

We then consider the generating function $$C(x,y,z,u)=\sum_{P\in \mathbf{C}}x^{\lth(P)}y^{\texttt{s}(P)}z^{\bar{\texttt{s}}(P)}u^{\last(P)}=\sum_{n,i,j,k\geq 0}c_{n,i,j,k}x^ny^iz^ju^k.$$ 

\begin{theorem} The generating function $C(x,y,z,u)$  satisfies
\begin{equation}\label{e1}
C(x,y,z,u)=xyu+\frac{xyu}{1-yu}\Bigl(C(x,z,y,1)-yu C(x,z,y,yu)\Bigr).
\end{equation}
\label{thm51}
\end{theorem}
\begin{proof}  If  $P$ has a single column, it consists of one cell, and its contribution to the generating function is $xyu$. Otherwise, assume that the polyomino has at least two columns. The polyomino $P$ is obtained from a polyomino $Q$ of length $\lth(Q)=\lth(P)-1$ by adding on its right a column with $\last(P)$ cells, by preserving the Catalan structure, that is,  the condition $\last(Q)+1\geq \last(P)$. Note that $\texttt{s}(P)=\bar{\texttt{s}}(Q)+\last(P)$ and  $\bar{\texttt{s}}(P)=\texttt{s}(Q)+\last(P)$.
    Then the contribution for these polyominoes is 
\begin{align*}
    &x\sum_{n,i,j,k} c_{n,i,j,k}\,x^{n}z^iy^j\left(uy+ u^2y^2+\cdots + u^{k+1}y^{k+1}\right)\\
    =&xuy\sum_{n,i,j,k} c_{n,i,j,k}\,x^{n}z^iy^j\left(\frac{1-u^{k+1}y^{k+1}}{1-uy}\right)\\
    =&\frac{xyu}{1-yu}\left(C(x,z,y,1)-yuC(x,z,y,yu)\right), 
\end{align*}
which completes the proof.
\end{proof}
By construction, the generating function $OD(x,y,1,1)$ obtained in the previous section (second coordinate of $\mathbf{V}(1)$ evaluated at $z=1$) coincides with the terms $x^{2n+1}y^i$ of $C(x,y,1,1)$. Similarly,   the generating function $EV(x,1,y,1)$ obtained in the previous section (first coordinate of $\mathbf{V}(1)$ evaluated at $y=1$) coincides with the terms $x^{2n}y^i$ of $C(x,y,1,1)$.
Thus, we have $C(x,y,1,1)=OD(x,y,1,1)+EV(x,1,y,1)$, and the result of the previous section provides a close form for the $C(x,y,1,1)$.


Below we show how we can obtain a closed form for the solution $C(x,y,z,u)$ of \eqref{e1}. Setting
\[
F(u)=C(x,y,z,u), \mbox{ and }  G(u)=C(x,z,y,u),
\]
Equation  \eqref{e1} becomes
\begin{equation}\label{e2}
F(u)=xyu+\frac{xyu}{1-yu}\Bigl(G(1)-yu G(yu)\Bigr).
\end{equation}
Interchanging $y$ and $z$ in \eqref{e1} gives the companion equation for $G$:
\begin{equation}\label{e3}
G(u)=xzu+\frac{xzu}{1-zu}\Bigl(F(1)-zu F(zu)\Bigr).
\end{equation}

Let $q=yz$. By \eqref{e3}, we have
\[
G(yu)=xqu+\frac{xqu}{1-qu}F(1)-\frac{xq^{2}u^{2}}{1-qu}F(qu).
\]
Multiplying this by $yu$ and substituting into \eqref{e2} yields the following lemma.

\begin{lemma}
The function $F(u)$ satisfies the equation
\begin{equation}\label{e10}
F(u)=A_y(u)+B_y(u) F(qu),
\end{equation}
where
\begin{align}
A_y(u)
&=xyu+\frac{xyu}{1-yu} G(1)
-\frac{x^2y^2 q u^3}{1-yu}
-\frac{x^2y^2 q u^3}{(1-yu)(1-qu)} F(1), \label{e5}\\
B_y(u)
&=\frac{x^2y^2 q^2 u^4}{(1-yu)(1-qu)}. \label{e6}
\end{align}
Similarly, $G$ satisfies
\[
G(u)=A_z(u)+B_z(u) G(qu),
\]
where $A_z(u),B_z(u)$ are obtained from \eqref{e5}--\eqref{e6} by interchanging $y$ and $z$ and $F(1)$ and $G(1)$.
\end{lemma}

\begin{lemma}
We have
\begin{equation}\label{e87}
F(u)=\sum_{n\ge 0}\left(\prod_{k=0}^{n-1}B_y(q^k u)\right) A_y(q^n u).
\end{equation}
\end{lemma}

\begin{proof}
Iterating \eqref{e10}, we inductively obtain, for every integer $N\geq 1$,
\[
F(u)=\sum_{n=0}^{N-1}\left(\prod_{k=0}^{n-1}B_y(q^k u)\right)A_y(q^n u)
+\left(\prod_{k=0}^{N-1}B_y(q^k u)\right)F(q^N u).
\] We claim that the second term vanishes as $N\to\infty$. Indeed, by \eqref{e6}, $B_y(u)$ has a factor $u^4$. Thus, $B_y(q^k u)$ has a factor $(q^k u)^4$ and therefore
$\prod_{k=0}^{N-1}B_y(q^k u)$ has a factor $u^{4N}$. Since $F(q^N u)$ is a formal power series in $u$, it has no negative powers of $u$. Thus,
\begin{equation}\label{ep}
\left(\prod_{k=0}^{N-1}B_y(q^k u)\right)F(q^N u)
\end{equation}
has a factor $u^{4N}$. It follows that for every nonnegative integer $m$, the  coefficient of $u^m$ in \eqref{ep} is $0$ for every $N>m/4$.
\end{proof}

Recall the $q$-Pochhammer symbol
\[
(a;q)_n=\prod_{k=0}^{n-1}(1-aq^k),
\]
where $n\geq 0$ and $(a;q)_0=1$.

\begin{lemma}
For $n\ge 0$,
\begin{equation}\label{e91}
\prod_{k=0}^{n-1}B_y(q^k u)
=
\frac{(x^2y^2u^4)^n q^{2n^2}}{(yu;q)_n (qu;q)_n}.
\end{equation}
\end{lemma}

\begin{proof}
By \eqref{e6},
\[
B_y(q^k u)
=\frac{x^2y^2u^4 q^{2+4k}}{(1-yq^k u)(1-q^{k+1}u)}.
\]
Thus
\[
\prod_{k=0}^{n-1}B_y(q^k u)
=
\frac{(x^2y^2u^4)^n
q^{\sum_{k=0}^{n-1}(2+4k)}
}{\prod_{k=0}^{n-1}(1-yq^k u)\;\prod_{k=0}^{n-1}(1-q^{k+1}u)}=\frac{(x^2y^2u^4)^n q^{2n^2}}{(yu;q)_n (qu;q)_n}.\qedhere
\]
\end{proof}

From \eqref{e5}, \eqref{e87}, and \eqref{e91}, we immediately deduce the following theorem.
\begin{theorem}
We have
\begin{equation}\label{ee1}
F(u)=\sum_{n\ge 0}\frac{(x^2y^2u^4)^n q^{2n^2}}{(yu;q)_n (qu;q)_n}
\left(
xy q^n u+\frac{xy q^n u}{1-yq^n u} G(1)
-\frac{x^2y^2 q^{1+3n}u^3}{1-yq^n u}
-\frac{x^2y^2 q^{1+3n}u^3}{(1-yq^n u)(1-q^{n+1}u)} F(1)
\right).
\end{equation}
The function $G(u)$ satisfies the analogous formula obtained by interchanging $y$ and $z$ and $F(1)$ and $G(1)$.
\end{theorem}

Let $P_y(u),Q_y(u)$, and $R_y(u)$, be functions such that $F(u)$ in \eqref{ee1} is written as
\begin{equation}\label{ee2}
F(u)=P_y(u)+Q_y(u) G(1)+R_y(u) F(1).
\end{equation}
Let $P_z(u),Q_z(u)$, and $R_z(u)$ be the corresponding functions for $G(u)$, i.e.,
\begin{equation}\label{ee22}
G(u)=P_z(u)+Q_z(u) F(1)+R_z(u) G(1).
\end{equation}
Substituting $u=1$ in \eqref{ee2} and \eqref{ee22} yields
\begin{align*}
F(1)&=P_y(1)+Q_y(1) G(1)+R_y(1) F(1),\\
G(1)&=P_z(1)+Q_z(1) F(1)+R_z(1) G(1).
\end{align*}
Thus, in matrix form,
\begin{equation}\label{er5}
\begin{pmatrix}
1-R_y(1) & -Q_y(1)\\
-Q_z(1) & 1-R_z(1)
\end{pmatrix}
\begin{pmatrix}F(1)\\ G(1)\end{pmatrix}
=
\begin{pmatrix}P_y(1)\\ P_z(1)\end{pmatrix}.
\end{equation}
If the determinant $\Delta$ of the matrix, which is given by
\[
\Delta=\bigl(1-R_y(1)\bigr)\bigl(1-R_z(1)\bigr)-Q_y(1)Q_z(1),
\]
is invertible, then
the system \eqref{er5} has the unique solution
\begin{equation}\label{ei9}
F(1)=\frac{P_y(1)\bigl(1-R_z(1)\bigr)+Q_y(1)P_z(1)}{\Delta},
\qquad
G(1)=\frac{P_z(1)\bigl(1-R_y(1)\bigr)+Q_z(1)P_y(1)}{\Delta}.
\end{equation}
Substituting \eqref{ei9} into \eqref{ee1} yields an explicit closed expression for
\[
C(x,y,z,u)=F(u)
\]

We obtained a closed form for the solution of the functional equation \eqref{e1}. We leave open the question of whether it is algebraic or $D$-finite, and if can obtain asymptotic approximation for the coefficient of $y^n$ in $C(1,y,1,1)$?

\section{Acknowledgement}
This research was funded, in part, by the Agence Nationale de la Recherche (ANR), grants ANR-22-CE48-0002, ANR-COMETA-GAE-25-CE48-0602 and by the Regional Council of Bourgogne-Franche-Comté.


\begin{thebibliography}{20}

\bibitem{Banfla}  C.~Banderier and P.~Flajolet. Basic analytic combinatorics of directed lattice paths. {\it Theoret. Comput. Sci.} \textbf{281} (2002), 37--80. \href{https://doi.org/10.1016/S0304-3975(02)00007-5}{https://doi.org/10.1016/S0304-3975(02)00007-5}.


\bibitem{Baril3} J.-L.~Baril, D.~Colmenares, J.~L.~ Ram\'irez, D.~Silva, L.~M.~Simbaqueba, and D.~Toquica. Consecutive pattern-avoidance in Catalan words according to the last symbol. {\it RAIRO Theor. Inform. Appl.} \textbf{58} (2024), Paper No. 1.  \href{https://doi.org/10.1051/ita/2024001}{https://doi.org/10.1051/ita/2024001}.

\bibitem{BGR} J.-L.~Baril, J.~F.~Gonz{\'a}lez, and J.~L.~Ram\'irez. Last symbol distribution in pattern avoiding Catalan words. {\it Math. Comput. Sci.} \textbf{18} (1) (2024).
\href{ https://doi.org/10.1007/s11786-023-00576-5}{ https://doi.org/10.1007/s11786-023-00576-5}.

\bibitem{Baril2} J.-L.~Baril, C.~Khalil, and V.~Vajnovszki. Catalan words avoiding pairs of length three patterns. {\it Discret. Math. Theor. Comput. Sci.} \textbf{22} (2) (2021), \# 5.  \href{https://doi.org/10.46298/dmtcs.6002}{ https://doi.org/10.46298/dmtcs.6002}


\bibitem{relation} J.-L.~Baril and J.~L.~Ram\'irez. Descent distribution on Catalan words avoiding ordered pairs of relations. {\it Adv. in Applied Math.} \textbf{149} (2023), 102551. \href{https://doi.org/10.1016/j.aam.2023.102551}{https://doi.org/10.1016/j.aam.2023.102551}
 

\bibitem{BLE3} A.~Blecher, C.~Brennan, and A.~Knopfmacher. Combinatorial parameters in bargraphs. {\it Quaest. Math.} \textbf{39} (2016), 619--635. \href{https://doi.org/10.2989/16073606.2015.1121932}{https://doi.org/10.2989/16073606.2015.1121932}


\bibitem{BlKn} A.~Blecher and A.~Knopfmacher. Cells of fixed height in Catalan words and restricted growth functions. {\it Adv. in Applied Math.}, \textbf{164} (2025), 102835. \href{https://doi.org/10.1016/j.aam.2024.102835}{
https://doi.org/10.1016/j.aam.2024.102835}

\bibitem{CallManRam} D.~Callan, T.~Mansour, and J.~L.~Ram\'irez.  Statistics on bargraphs of Catalan words. {\it J. Autom. Lang. Comb.} \textbf{26} (2021), 177--196.
\href{https://doi.org/10.25596/jalc-2021-177}{
https://doi.org/10.25596/jalc-2021-177}

\bibitem{Chen} S.~E.~Cheng, D.~P.~Eu, and T.~S.~Fu. Area of Catalan paths on a checkerboard. {\it European J Combin.} \textbf{28}(4) (2007), 1331--1344. \href{https://doi.org/10.1016/j.ejc.2006.01.006}{https://doi.org/10.1016/j.ejc.2006.01.006}


\bibitem{DCastro} R.~De~Castro, A.~Ramírez, and J.~L.~Ramírez. Applications in enumerative combinatorics of infinite weighted automata and graphs. \emph{Sci. Ann. Comput. Sci.} Vol. XXIV (1), (2014), 137--171, 
\href{https://doi.org/10.7561/SACS.2014.1.137}{
https://doi.org/10.7561/SACS.2014.1.137}

\bibitem{Deu} E.~Deutsch. Dyck path enumeration. {\it Discrete Math.} \textbf{204} (1999), 167--202. \href{https://doi.org/10.1016/S0012-365X(98)00371-9}{https://doi.org/10.1016/S0012-365X(98)00371-9}


\bibitem{fried} S.~Fried. Black-white cell capacity in $k$-ary words and permutations. \href{https://arxiv.org/abs/2509.07533}{https://arxiv.org/abs/2509.07533}.

\bibitem{gasper} G.~Gasper and M.~Rahman. \emph{Basic Hypergeometric Series}. Cambridge University Press, 1990.
 
 \bibitem{Book1} A.~J.~Guttmann (Ed.) \emph{Polygons, Polyominoes and Polycubes}.  Lecture Notes in Physics 775. Springer, Heidelberg, Germany, 2009.
 

\bibitem{ManSha2} T.~Mansour and A.~Sh.~Shabani. Enumerations on bargraphs. {\it Discrete Math. Lett.} \textbf{2} (2019), 65--94.

\bibitem{ManRam} T.~Mansour and J.~L.~Ram\'irez. Enumerations on polyominoes determined by Fuss-Catalan words. {\it Australas. J. Comb.} \textbf{81} (2021), 447--457.

\bibitem{Toc}
 T.~Mansour,  J.~L.~Ram\'irez, and D.~A. Toquica. Counting lattice   points on bargraphs of {C}atalan words.  {\it Math. Comput. Sci.} \textbf{15} (2021), 701--713. \href{https://doi.org/10.1007/s11786-021-00501-8}{https://doi.org/10.1007/s11786-021-00501-8}

\bibitem{OEIS} OEIS Foundation Inc., The On-Line Encyclopedia of Integer Sequences, \url{http://oeis.org/}.
 
\bibitem{Prod} H.~Prodinger. Motzkin paths of bounded height with
two forbidden contiguous subwords of length two. J. Integer Seq. \textbf{28} (2025), Article 25.7.6. 

\bibitem{RaissouliKacha2000} M.~Raissouli and A.~Kacha. Convergence of matrix continued fractions. {\it Linear Algebra Appl.}  \textbf{120} (1-3) (2000), 115--129. \href{https://doi.org/10.1016/S0024-3795(00)00196-8}{https://doi.org/10.1016/S0024-3795(00)00196-8}

\bibitem{AlejaRam} J.~L.~Ram\'irez and A.~Rojas-Osorio. Consecutive patterns in Catalan words and the descent distribution. {\it Bol. Soc. Mat. Mex.} \textbf{29} (2023), Article \#60.  \href{https://doi.org/10.1007/s40590-023-00532-0}{https://doi.org/10.1007/s40590-023-00532-0}.
 
\bibitem{sorokin} V.~N.~Sorokin. Matrix continued fractions. {\it J. Approx. Theory} \textbf{96} (1999), 237--257.
\href{https://doi.org/10.1006/jath.1998.3232}{https://doi.org/10.1006/jath.1998.3232}

\bibitem{stan} R.~P.~Stanley. \emph{Catalan Numbers}.  Cambridge University Press, 2015.


\end{thebibliography}
\end{document}